\documentclass[reqno]{amsart}
\usepackage{amssymb,amsthm,amsfonts,amstext}
\usepackage{amsmath}
\usepackage[notcite,notref]{showkeys}
\usepackage{zito}

\makeatletter \@addtoreset{equation}{section} \makeatother

\renewcommand\thefigure{\thesection.\@arabic\c@figure}
\renewcommand\thetable{\thesection.\@arabic\c@table}
\newtheorem{theorem}{Theorem}[section]
\newtheorem{lemma}[theorem]{Lemma}
\newtheorem{proposition}[theorem]{Proposition}

\newtheorem{remark}[theorem]{Remark}

\newcommand{\mc}[1]{{\mathcal #1}}

\newcommand{\bb}[1]{{\mathbb #1}}

\newcommand{\<}{\langle}
\renewcommand{\>}{\rangle}

\title{Hydrodynamic limit of particle systems with long jumps}
\date{}
\author{Milton Jara}
\address{CEREMADE, Universit\'e Paris-Dauphine\\
Place du Mar\'echal de Lattre De Tassigny, Paris CEDEX 75775, France
}

\begin{document}

\begin{abstract}

We consider some interacting particle processes with long-range dynamics: the zero-range and exclusion processes
with long jumps. We prove that the hydrodynamic limit of these processes corresponds to a (non-linear in general) fractional heat equation. The scaling in this case is superdiffusive. In addition, we discuss a central limit theorem for a tagged particle on the zero-range process and existence and uniqueness of solutions of the Cauchy problem for the fractional heat equation.

\end{abstract}


\renewcommand{\subjclassname}{\textup{2000} Mathematics Subject Classification}

\address{Universit\'e Catholique de Louvain, 2 Chemin du Cyclotron, B-1348, Louvain-la-Neuve, Belgium} 
\email{jara@ceremade.dauphine.fr}

\maketitle

\section{Introduction}

Interacting particle systems have been the subject of intense study for the last 40 years. This is due to the fact that, in one hand, they present many of the collective features found in real physical systems, and in the other hand they are, up to some extent, mathematically tractable. The rigorous study of interacting particle systems has lead in many cases to a more detailed understanding of the microscopic mechanism behind those collective phenomena. We refer to \cite{KL} for further references and to \cite{BDGJL} for a recent result which we think is a good example of the success   this plan. 

Since their introduction by Spitzer, the zero-range process and the exclusion process have been among the most studied interacting particle systems, and they have served as a test field for new mathematical and physical ideas.

During the last years, and specially due to applications in finance, fractional Laplacians $L =-(-\Delta)^s$, $s \in (0,1)$ and their probabilistic counterparts, the L\'evy processes, have received grown attention. Some of the key properties of L\'evy processes are the presence of jumps, the lack (at least in the cases associated to $L$) of bounded variance and a super-diffusive behavior. 

From the point of view of statistical mechanics, it is desirable to have a derivation of the partial differential equations ruling the evolution of these super-diffusive systems, which usually involve the fractional Laplacian, from microscopic models. We are not aware of any work on this direction, so we have decided to provide in this article such a derivation. Our aim is to obtain the hydrodynamic limit of interacting particle systems on which particle may perform long jumps, in the context of the exclusion and zero-range processes.  For these models, the corresponding hydrodynamic equation is given by a fractional (non linear in general) heat equation of the form $\partial_t u = L u$, where $L$ is the generator of a symmetric, $\alpha$-stable L\'evy process and in particular it includes the fractional Laplacians $-(-\Delta)^s$ as particular cases. We believe that this equation will emerge as well as the hydrodynamic limit of other particle systems which are superdiffusive in nature, like heat conduction models with conservative noise in dimension $d=1$, Ornstein-Uhlenbeck particles in laminar flow, etc.

This article is organized as follows. In Section \ref{s1} we give detailed definitions of the models considered here. In Section \ref{s2} we prove our main result for the exclusion process with long jumps. At the current state of art, this result is more or less elementary. We included a complete proof for two reasons. In general terms, the theory of hydrodynamic limits of particle systems is a hard subject for the non-expert reader, since it involves a mixture of purely probabilistic ideas and ideas coming directly from analysis and PDE theory. Therefore, we believe that it is a good way to introduce the subject in a more elementary way, simplifying later the exposition for the zero-range process. And secondly, the results of this section were used without proof in \cite{Jar} to obtain an invariance principle for a tagged particle in the exclusion process with long jumps. In Section \ref{s3} we prove the main theorem for the zero-range process, leaving various technical results for the following sections. In Section \ref{s4} we prove the key technical input, known in the literature as the {\em replacement lemma}. In Section \ref{s5} we prove the so-called {\em moving particle lemma}, which we believe is the main feature that differentiates super-diffusive systems from the most studied diffusive systems, and the new proof of this lemma is the main technical  novelty of this work. Although a portion of the exposition is by now classical, the tools needed have been gathered from many different places, and we have decided to include detailed proofs of most of the propositions taken from elsewhere in order to keep the exposition as clear as possible. We invite the most specialized reader to skip the more standard parts of the exposition. 

In Section \ref{s6} we prove the energy estimate, which is crucial in order to obtain a uniqueness criterion for the hydrodynamic equation. In this section we introduce a variational formula for the Fisher information for the fractional heat equation which seems to be new in the literature. This formula involves a natural generalization of the space of test functions to antisymmetric functions of two space coordinates. In the derivation of hydrodynamic limits of particle systems, a key analytical input is an existence result for the corresponding hydrodynamic equation. After consulting some experts in the field, it seems that the required uniqueness results are not available in the literature, so in Section \ref{s7} we obtain some uniqueness results tailored to our needs. Those results may be of independent interest, and they are independent of the rest of paper. The proof in the linear case is due to Luis Silvestre. 

In Section \ref{s8} we prove an invariance principle for a tagged particle in the one-dimensional zero-range process, following recent results in \cite{Jar}, \cite{JLS}. The limiting process is a time-inhomogeneous process of independent increments, related to the solution of the hydrodynamic equation. In Appendix \ref{A2} Êwe explain how to deal with the borderline case $\alpha =2$, which leads to the usual heat equation, but with a superdiffusive time scaling. The hydrodynamic limit of the zero-range process in its full generality is proved only under a restrictive attractiveness condition for the system. However, we point out that it is only at the level of the uniqueness criterion for the hydrodynamic equation that this result is needed. Certainly most of the proofs can be obtained in a simpler way using coupling arguments, available only for attractive systems. In this article, attractiveness is only used to obtain a bound for the energy of the solutions in terms of their Fisher information. 
We conjecture that, like in the diffusive case, uniqueness under the weaker bound on Fisher information holds. In Appendix \ref{B} we explain how to consider unbounded initial profiles in the formulation of the hydrodynamic limit for attractive systems. In Appendix \ref{C} we explain how to handle general bounded initial profiles. Both results are based on coupling tecniques which are more or less standard in the literature. A combination of both appendices can be used to handle initial conditions of the form $u_0(x)+v_0(x)$, where $u_0$ is bounded and $v_0$ is in $\mc L^1(\bb R^d)$. 

\section{The models}
\label{s1}
Consider the integer lattice $\bb Z^d$, $d \geq 1$ and let $p: \bb Z^d \to [0,\infty)$ be such that $\sum_{z \in \bb Z^d} p(z) =p^*<+\infty$. We call $p(\cdot)$ the {\em transition rate}. The exclusion process associated to $p(\cdot)$ is the Markov process $\eta_t$ defined in $\Omega_{ex} = \{0,1\}^{\bb Z^d}$ and generated by the operator
\[
L^{ex} f(\eta) = \sum_{x,y \in \bb Z^d} p(y-x) \eta(x) \big(1-\eta(y)\big) \big[f(\eta^{xy}) - f(\eta)\big].
\]

Here $\eta$ denotes a generic element of $\Omega_{ex}$, $f: \Omega_{ex} \to \bb R$ is a local function, that is, it depends on the value of $\eta(z)$ for only a finite number of points $z \in \bb Z^d$, and
\[
\eta^{xy} (z) = 
\begin{cases}
\eta(y), & z=x\\
\eta(y), & z=y\\
\eta(z), &z \neq x,y.\\
\end{cases}
\]

For details about the construction and properties of the process $\eta_t$, we refer to \cite{Lig}. The dynamics of this process is easy to describe. Initially, particles are distributed in $\bb Z^d$ in such a way that there is at most one particle per site. Each particle, independent of the other particles, waits an exponential time of rate $p^*$, at the end of which it picks a site $y \in \bb Z^d$ with probability $p(y-x)/p^*$, where $x$ is the current position of the particle. If the site $y$ is empty, the particle jumps from $x$ to $y$. Otherwise the particle stays at $x$. In any case, a new exponential time starts afresh. Notice that for initial configurations with a finite number of particles, the process $\eta_t$ is just a system of independent random walks with transition rate $p(\cdot)$, conditioned to have at most one particle per site. When the number of particles is infinite, the construction of the process has been carried out by Liggett \cite{Lig2}.

Denote by $\bb N_0$ the set of non-negative integers and let $g: \bb N_0 \to [0,\infty)$ be a function such that $g(0)=0$. We assume that $g(n)>0$ for $n>0$. The zero-range process with interaction rate $g(\cdot)$ and transition rate $p(\cdot)$ is the Markov process $\xi_t$ defined in $\Omega_{zr}^0 = (\bb N_0)^{\bb Z^d}$ and genetared by the operator
\[
L^{zr} f(\xi) = \sum_{x,y \in \bb Z^d} p(y-x) g\big(\xi(x)\big) \big[f(\xi^{x,y}) -f(\xi)\big].
\]

Here $\xi$ is a generic element of $\Omega_{zr}^0$, $f: \Omega_{zr}^0 \to \bb R$ is a ``suitable'' function and
\[
\xi^{x,y}(z) = 
\begin{cases}
\xi(x)-1, & z=x\\
\xi(y)+1, &z=y\\
\xi(z), &z \neq x,y.\\
\end{cases}
\]

The exact meaning of  ``suitable'' depends on the choice of $p(\cdot)$, but for local functions $f$, a Lipschitz condition of the form
\[
\big|f(\xi^{x,y})-f(\xi)\big| \leq C(f)
\]
for any $\xi,x,y$ will be sufficient.

The dynamics of this process is the following: at a site $x \in \bb Z^d$ and independently for each site, the particles wait an exponential time of rate $p^*g(\xi(x))$, where $\xi(x)$ denotes the number of particles at site $x$. At the end of this exponential time, one of the particles at $x$ jumps to a site $y$, randomly chosen with probability $p(y-x)/p^*$. Then, a new exponential time starts afresh. Notice that the particles interact between them only when they are at the same site, explaining the denomination ``zero-range process''. The case $g(n)=n$ corresponds to a system of independent random walks in $\bb Z^d$ with transition rate $p(\cdot)$.

Differently from the exclusion process, the zero-range process can present explosions if the initial configuration of particles has too many particles at infinity.  In this article, we assume that the interaction rate satisfies the Lipschitz  condition stated above: there is a finite constant $\kappa>0$ such that $|g(n+1)-g(n)| \leq \kappa$ for any $n \geq 0$. We refer to this condition as {\bf (LG)}. Under {\bf (LG)}, the process $\xi_t$ is well-defined for any bounded initial condition and also $a.s.$ with respect to any initial measure $\mu$ in $\Omega_{zr}^0$ such that $\sup_x E_\mu[\xi(x)]<+\infty$. 
Fore a more detailed discussion about this point, the construction of the process and related topics, we refer to \cite{And}.

\subsection{Invariant measures}
\label{s1.1}
Let $\rho$ be a fixed number in $[0,1]$. Define $\mu_\rho$ as the product measure in $\Omega_{ex}$ with one-site marginals satisfying
\[
\mu_\rho(\eta(x)=1) = 1 - \mu_\rho(\eta(x)=0) = \rho.
\]

It is well known \cite{Lig} that the measures $\{\mu_\rho; \rho \in [0,1]\}$ are invariant under the evolution of $\eta_t$. When $p(\cdot)$ is symmetric, the measures $\{\mu_\rho\}_\rho$ are also reversible. And if the transition rate $p(\cdot)$ is irreducible, then the measures $\{\mu_\rho\}_\rho$ are also ergodic. Notice that $\int \eta(x) d\mu_\rho = \rho$, that is, the density of particles per site is equal to $\rho$ for the measure $\mu_\rho$. The fact that there exists a family $\{\mu_\rho\}_\rho$ of invariant measures parametrized by the density o particles reflects the fact that the dynamics is conservative: particles are neither created nor destroyed by the dynamics.

For the zero-range process there also exists a family of invariant measures of product form \cite{And}. Take $\phi \geq 0$ and define $\bar \nu_\phi$ as the measure in $\Omega_{zr}^0$ of marginals given by
\[
\bar \nu_\phi \big(\xi(x_1)=k_1,\dots,\xi(x_l)=k_l\big)
	= \prod_{i=1}^l \frac{1}{Z(\phi)} \frac{\phi^{k_i}}{g(k_i)!},
\]
where $g(k)!= g(1)\cdots g(k)$ for $k>0$, $g(0)! = 1$ and $Z(\phi)$ is the normalization constant
\[
Z(\phi) = \sum_{k \geq 0} \frac{\phi^k}{g(k)!}.
\]

Notice that $Z(\phi)$ is an increasing function of $\phi$. Therefore, there is a maximal value $\phi_c$ (perhaps equal to $+\infty$) such that $Z(\phi)$ is finite for $\phi<\phi_c$ and $Z(\phi)$ is infinite for $\phi > \phi_c$. If there is a positive constant $\epsilon_0$  such that $g(n) \geq \epsilon_0$ for any $n>0$, $\phi_c>0$ and $Z(\phi_c)=+\infty$. The family of measures $\{\bar \nu_\phi; \phi <\phi_c\}$ is invariant under the dynamics of $\xi_t$. When the transition probability $p(\cdot)$ is symmetric, the measures $\bar \nu_\phi$ are reversible. And when $p(\cdot)$ is irreducible, the measures $\bar \nu_\phi$ are ergodic. 

The zero-range process also conserves the number of particles, so it would be more natural to parametrize the invariant measures by the density of particles per site. Define $\rho(\phi) = \int \xi(x) \bar \nu_\rho$. It is not difficult to see that $\rho(\phi)$ is a differentiable, strictly increasing function from $[0,\phi_c)$ to $[0,\infty)$. Let us write $\rho_c = \lim_{\varphi  \to \varphi_c} \rho(\varphi)$. The inverse function $\phi(\rho)$ is well defined for any $\rho \in [0,\rho_c)$. We define $\nu_\rho = \bar \nu_{\phi(\rho)}$. Notice that $\phi(\rho) = \int g(\xi(x)) \nu_\rho$.

The occupation variables $\xi(x)$ and the interaction rates $g(\xi(x))$, $x \in \bb Z^d$ have exponential moments of sufficiently small order. In fact, for $\theta \in \bb R$,
\[
M_\rho(\theta) =\int e^{\theta \xi(x)} d\bar \nu_\phi = \frac{1}{Z(\phi)} \sum_{k \geq 0} \frac{(\phi e^\theta)^k}{g(k)!}
	= \frac{Z(\phi e^\theta)}{Z(\phi)}
\]
and we conclude that $M_\rho(\theta)$ is finite for $\theta <\log(\phi_c/\phi)$. From {\bf (LG)}, we obtain the bound $g(n) \leq \kappa n$, and therefore $\int e^{\theta g(\xi(x))} d\bar \nu_\phi$ is finite if $\theta \leq \kappa^{-1} \log(\phi_c/\phi)$. We say that the interaction rate $g(\cdot)$ satisfies {\bf (FEM)} if $M_\rho(\theta)$ is finite for every $\theta, \rho\geq0$. A simple application of the ratio test shows that a non-decreasing interaction rate $g(\cdot)$ satisfies {\bf (FEM)} if and only if $\phi_c =+\infty$, and equivalently if and only if $\lim_{n \to \infty} g(n) =+\infty$ (this last limit always exists, since $g(\cdot)$ is non-decreasing). We say that $g(\cdot)$ satisfies {\bf (B)} if $g(\cdot)$ is non-decreasing and if it does not satisfy {\bf (FEM)}. Therefore, $g(\cdot)$ satisfies {\bf (B)} if and only if $g(\cdot)$ is non-decreasing and bounded.

We already mentioned that the process $\xi_t$ is well defined $a.s.$ with respect to an initial measure $\mu$ such that $\sup_x E_\mu[\xi(x)]<+\infty$. Here we quote a more precise statement, due to Andjel \cite{And}:

\begin{proposition}
\label{p1} For any transition rate $p(\cdot)$, there exists a function $\sigma: \bb Z^d \to [0,\infty)$ such that:
\begin{itemize}
\item[i)] $\sum_z \sigma(z) <+\infty$,
\item[ii)] The zero-range process $\xi_t$ is a strong Markov process when defined on the set
\[
\Omega_{zr} = \big\{ \xi \in \Omega_{zr}^0 ; \sum_{z \in \bb Z^d} \xi(z) \sigma(z) <+\infty\big\}.
\]
\end{itemize}
\end{proposition} 

Notice that $\Omega_{zr}$ has full measure for any of the invariant measures $\nu_\rho$, $\rho \geq 0$.
From now on, we will always define the process $\xi_t$ in $\Omega_{zr}$.

\begin{remark}
For $x \in \bb R^d$, we denote by $||x||$ the Euclidean norm $(x_1^2+\dots+x_d^2)^{1/2}$ of $x$, and by $|x|$ the supremum norm $\sup\{|x_1|,\dots,|x_d|\}$. For a given Polish space $\mc E$, we denote by $\mc D([0,\infty),\mc E)$ ($\mc D([0,T],\mc E)$ resp.) the space of c\`adl\`ag functions $f: [0,\infty) \to \mc E$ ($f: [0,T] \to \mc E$ resp.) equipped with the $J$-Skorohod topology.
\end{remark}

\subsection{Homogeneous transition rates}
When the transition rate $p(\cdot)$ has mean zero and finite range, that is, when  $\sum_z zp(z) =0$ and $p(z) =0$ for $z$ big enough, it is well known that the hydrodynamic limit of the processes $\eta_t$, $\xi_t$ is diffusive and given by the heat equation in the case of the exclusion process $\eta_t$ and by a nonlinear heat equation for the process $\xi_t$ (see \cite{KL} and the references therein). In the case $\sum_z z p(z) =m \neq 0$, the hydrodynamic limit appears in the hyperbolic scale and corresponds to a conservation law. Therefore, in order to obtain a superdiffusive scaling limit, it is necessary for the transition rate $p(\cdot)$ to have arbitrarily long jumps. Probably, the most natural choice for such a $p(\cdot)$ should be to take $p(x) = 1/ ||x||^{d+\alpha}$, where $||x||=(x_1^2+ \dots+ x_d^2)^{1/2}$ is the Euclidean norm of $x$ and $\alpha>0$ is an arbitrary constant. Notice that $p^*<+\infty$ due to the condition $\alpha>0$. Now we describe a broader class of transition rates for which $1/||x||^{d+\alpha}$ will be a particular case. A function $h: \bb R^d \setminus \{0\} \to \bb R$ is said to be {\em homogeneous} of degree $\beta \in \bb R$ if for any $x \in \bb R^d \setminus \{0\}$ and any $\lambda >0$, $h(\lambda x) = \lambda^\beta h(x)$. We say that a transition rate $p(\cdot)$ is homogeneous, regular of degree $\beta$ if there is a function $h: \bb R^d\setminus \{0\} \to \bb R$ homogeneous of degree $\beta$, continuous and strictly positive such that $p(x) = h(x)$ for any $x \in \bb Z^d \setminus \{0\}$.

We assume the following hypothesis on $p(\cdot)$:

\begin{description}
\item[(P)] The transition rate $p(\cdot)$ is symmetric and homogeneous of degree $-(d+\alpha)$, with $\alpha \in (0,2)$.
\end{description}

The restriction $\alpha >0$ comes from the fact that $\sum_z p(z) <+\infty$, and in particular the function $h(\cdot)$ has to be integrable outside a ball around the origin. The restriction $\alpha <2$ comes from the fact that for $\alpha >2$, $\sum_z z^2 p(z) <+\infty$ and in that case the hydrodynamic limits of $\eta_t$ and $\xi_t$ are still diffusive. The boundary case $\alpha =2$ is special, since the hydrodynamic limit is expected to be the usual heat equation, but the scaling contains a logarithmic correction (see Appendix \ref{A2}).

Notice that in $d=1$, the unique transition rates $p(\cdot)$ satisfying {\bf (P)} are $p(z) = c/||z||^{1+\alpha}$, $c>0$. In fact, the class of transition rates $p(\cdot)$ satisfying {\bf (P)} is homeomorphic to the class of measures $m$ in the sphere $\bb S^{d-1}=\{x \in \bb R^d; ||x||=1\}$ with a continuous, strictly positive density with respect to Lebesgue measure and satisfying $m(A)=m(-A)$ for any Borel set $A \in \bb S^{d-1}$.
We call the processes $\eta_t$  and $\xi_t$ associated to a transition rate $p(\cdot)$ satisfying {\bf (P)} the exclusion process and zero-range process with long jumps.

\subsection{The hydrodynamic limit for $\eta_t$}
Let $p(\cdot)$ be a transition rate satisfying {\bf (P)} and let $h(\cdot)$, $-(d+\alpha)$ be the corresponding homogeneous function and degree. These parameters will be fixed throughout the rest of this article. We start defining the pseudo-differential operator
\[
\mc L F(x) = \frac{1}{2} \int_{\bb R^d} h(y) \big(F(x+y)+F(x-y)-2F(x)\big) dy
\]
for $F \in \mc C_c^2(\bb R^d)$, the set of twice-continuously differentiable functions of bounded support. Since the function $F$ is bounded, the integral is absolutely convergent outside a ball around the origin (at this point we need $\alpha>0$). And using a second-order Taylor expansion of $F$ around $x$ we see that the integral is also absolutely convergent around 0 (here we need $\alpha <2$). In fact, $\mc L F: \bb R^d \to \bb R^d$ is a continuous function. Notice that in the case $h(y) = c/||y||^{d+\alpha}$, the operator $\mc L$ corresponds to a constant multiple of the fractional Laplacian $-(-\Delta)^{\alpha/2}$.

Let $u_0: \bb R^d \to \bb R$ be a measurable function. We say that a sequence of probability measures $\{\mu^n\}_n$ in $\Omega_{ex}$ is {\em associated} to the initial profile $u_0$ if for any function $F$ in the set $\mc C_c(\bb R^d)$ of continuous functions of bounded support and every $\epsilon >0$,
\[
\lim_{n \to \infty} \mu^n\Big\{ \eta; \Big| n^{-d} \sum_{z \in \bb Z^d} \eta(z) F(z/n) - \int u_0(x) F(x) dx\Big| >\epsilon \Big\} =0.
\]

We adopt the same definition for a sequence of measures $\{\nu^n\}_n$ in $\Omega_{zr}$, exchanging $\eta$ by $\xi$ in the previous relation. Notice that in order to have a sequence of measures $\{\mu^n\}_n$ associated to an initial profile $u_0$, it is necessary to have $0 \leq u_0(x) \leq 1$ for any $x \in \bb R^d$.

Fix $T >0$. We say that a measurable function $u:[0,T] \times \bb{R}^d \to \bb R$ is a {\em weak solution} of the Cauchy problem
\begin{equation}
 \label{ec1}
\left\{
\begin{array}{rcl}
\partial_t u_t &=& \mc L u_t\\
u(0,\cdot) &=& u_0(\cdot)\\
\end{array}
\right.
\end{equation}
if $||u||_{\infty,T} = \sup\{u(t,x); t \in [0,T], x \in \bb R^d\}$ is finite and for any smooth function $G: [0,\infty) \times \bb{R}^d \to \bb R$ of compact support we have 
\[
 \int_{\bb R^d} u(T,x) G(T,x) dx -\int_{\bb R^d} u_0(x) G(0,x) dx 
 	- \int_0^T\int_{\bb R^d}u(t,x) \big\{\partial_t + \mc L\big\}G(t,x) dx dt =0.
\]

At this point we are ready to define what we mean by the hydrodynamic limit of the process $\eta_t$.

\begin{theorem}
\label{t1}
Let $u_0: \bb R^d \to [0,1]$ be an initial profile, and let $\{\mu^n\}$ be associated to $u_0$. Define $\eta_t^n$ as the speeded-up process $\eta_{t n^\alpha}$ starting from the initial measure $\mu^n$ and let $\mu^n(t)$ be the distribution in $\Omega_{ex}$ of $\eta_t^n$. Then, the sequence $\{\mu^n(t)\}_n$ is associated to the function $u(t,\cdot)$, where $u(t,x)$ is the unique weak solution of the Cauchy problem (\ref{ec1}).
\end{theorem}

Equation (\ref{ec1}) is known as the {\em hydrodynamic limit} of the process $\eta_t^n$. It is important to notice the superdiffusive scaling $n^\alpha$ in this theorem. For a given measurable initial profile $u_0: \bb R^d \to [0,1]$, it is not difficult to construct a sequence of measures $\{\mu^n\}_n$ associated to it. In fact, it is enough to consider the product measures $\mu^n$ defined by
\begin{equation}
\label{ec1.5}
 \mu^n \big\{ \eta(z)=1\big\} = 1-\mu^n \big\{ \eta(z)=0\big\} = u_0^n(z)=:\int\limits_{|x-\frac{z}{n}|\leq\frac{1}{2n}} u_0(x) dx.
\end{equation}

Let $\mc M_+(\bb R^d)$ be the space of non-negative, Radon measures in $\bb R^d$. 
Let the {\em empirical density} be the random measure in $\mc M_+(\bb R^d)$ defined by
\[
 \pi^n(dx) = n^{-d} \sum_{z \in \bb Z^d} \eta(z) \delta_{z/n}(dx),
\]
where $\delta_x$ is the Dirac-$\delta$ distribution at $x$. It is not difficult to see that the sequence $\{\mu^n\}_n$ is associated to $u_0$ if and only if the sequence $\{\pi^n(dx)\}_n$ of random measures satisfies a weak law of large numbers with respect to the weak topology in $\mc M_+(\bb R^d)$, with limit measure $u_0(x)dx$. 

\subsection{Atractiveness and stochastic domination}
\label{s1.4}
Let $\xi$, $\xi'$ be two elements of $\Omega_{zr}$. We say that $\xi \preceq \xi'$ if and only if $\xi(x) \leq \xi'(x)$ for any $x \in \bb Z^d$. The relation $\preceq$ defines a partial order in $\Omega_{zr}$.  

We say the the process $\eta_t$ is {\em attractive} if for any two initial configurations of particles $\xi^1$, $\xi^2$ in $\Omega_{zr}$ with $\xi^1 \preceq \xi^2$, there exists a process $(\xi^1_t, \xi^2_t)$ in $\Omega_{zr} \times \Omega_{zr}$ such that the evolution of $\xi^i_t$, $i=1,2$ corresponds to a zero-range process starting from $\xi^i$ and $\xi_t^1 \preceq \xi^2_t$ $a.s.$ for any $t \geq 0$.

It is not hard to see that the process $\xi_t$ is atractive when the interaction rate $g(\cdot)$ is non-decreasing. The construction of the process $(\xi^1_t,\xi^2_t)$ is simple. First we take a zero-range process $\xi^1_t$ starting from $\xi^1$. Fix now a realization of the process $\xi^1_t$ and let $\delta_t$ be the zero-range process with initial configuration $\delta$ defined by $\delta(x) = \xi^2(x) - \xi^1(x)$ and inhomogeneous interaction rate $g_{x,t}(\delta) = g(\delta(x)+\xi^1_t(x)) -g(\xi^1_t(x))$. The process $\delta_t$ is well defined exactly due to the fact that $g(\cdot)$ is non-decreasing. Then we define $\xi^2_t$ by taking $\xi^2_t(x) = \xi_t^1(x)+\delta_t(x)$. For a justification of the validity of this construction, see Section 4 of \cite{And}.

Let $\nu$, $\nu'$ be two probability measures in $\Omega_{zr}$. We say that $\nu$ is {\em  stochastically dominated} by $\nu'$ and we write $\nu \preceq \nu'$, if there exists a measure $\hat \nu$ in $\Omega_{zr} \times \Omega_{zr}$ such that $\nu(\cdot) = \hat \nu(\cdot, \Omega_{zr})$, $\nu'(\cdot) = \hat \nu(\Omega_{zr},\cdot)$ and $\hat \nu (\xi^1 \preceq \xi^2) =1$.

We say that a function $F: \Omega_{zr} \to \bb R$ is non-decreasing if $F(\xi) \leq F(\xi')$ whenever $\xi \preceq \xi'$. Taking the expectation of $F(\xi^1) - F(\xi^2)$ with respect to $\hat \nu$, it is easy to see that if $\nu \preceq \nu'$, then 
\[
\int F d\nu \leq \int F d\nu'
\]
for any non-decreasing function $F$. In fact, the validity of this last relation for any non-decreasing, bounded function $F$ is sometimes used as the definition of stochastic domination.

The following property of the process $\xi_t$ follows easily from its atractiveness.

\begin{proposition}
\label{p2}
Let $\nu^1$, $\nu^2$ be two probability measures defined in $\Omega_{zr}$ with $\nu^1 \preceq \nu^2$. Denote by $\nu^i(t)$ the distribution at time $t$ of the zero-range process $\xi_t$ starting from the initial measure $\nu^i$, $i=1,2$. Then $\nu^1(t)$ is stochastically dominated by $\nu^2(t)$. 
\end{proposition}

For any $\rho^1\leq \rho^2$, we have $\nu_{\rho^1} \preceq \nu_{\rho^2}$ (see Section 1 of \cite{KL}). Here we give  a quick probabilistic proof of this fact. It is enough to prove that $q_{\rho^1} \preceq q_{\rho^2}$, where $q_\rho$ is the probability measure in $\bb N_0$ given by the one-site marginal of $\nu_\rho$. Notice that $q_\rho$ is the unique invariant measure of the {\em ladder process} in $\bb N_0$ defined as follows. A particle at site $n \geq 0$ goes down with exponential rate $g(n)$ and goes up with exponential rate $\phi(\rho)$. Now take two particles $X^1$, $X^2$ evolving in $\bb N_0$ as follows. At time $t = 0$ they both start at $n = 0$. When they are not together they evolve independently, following a ladder process of rate $\phi(\rho^i)$, $i=1, 2$. When they are together, say at site $n$, they wait an exponential time of rate $g(n)$ at the end of which they go one step down together. They also wait an exponential time of rate $\phi(\rho^2)$, at the end of which particle $X^2$ goes one step up, and particle $X^1$ goes  up together with particle $X^2$ with probability $\phi(\rho^1)/\phi(\rho^2)$ and it stays at $n$ with probability $1 - \phi(\rho^1)/\phi(\rho^2)$. It is clear that the law of $X^i$ corresponds to a ladder process of rate $\phi(\rho^i)$, $i=1,2$ and $X^1(t) \leq X^2(t)$ for any $t \geq 0$. The process $(X^1(t),X^2(t))$ constructed in this way is recurrent. In fact, $X^2(t)$ is recurrent, so the return time to 0 is $a.s.$ finite and has finite mean. But each time $X^2(t)$ returns to $0$, $X^1(t)$ also returns to $0$. Therefore, the return time of $(X^1(t),X^2(t))$ to $(0,0)$ is $a.s.$ finite and has finite mean. We conclude that $(X^1(t),X^2(t))$ is recurrent, and therefore it has a unique invariant measure. For a process with only one invariant measure, the C\`esaro means of the distributions at time $t$ of the process converge to the invariant measure. In our case, the C\`esaro means of the distributions of $(X^1(t),X^2(t))$ converge to a probability measure $\hat q$ with marginals $q_{\rho^i}$ and satisfying $\hat q(x^1 \leq x^2) =1$.

\subsection{Hydrodynamic limit for $\xi_t$}
\label{s1.5}
Before stating our result about the hydrodynamic limit of the zero-range process $\xi_t$, we need some definitions. Let $\mu$, $\nu$ be two probability measures defined in some  measurable space $\Omega$. The relative entropy of $\mu$ with respect to $\nu$ is defined by
\[
H(\mu|\nu) = \int \frac{d\mu}{d\nu} \log \frac{d\mu}{d\nu} d\nu
\]
if $\mu$ is absolutely continuous with respect to $\nu$, and $H(\mu|\nu) = +\infty$ otherwise.

Remember the definition of $\phi(\rho)$ in Section \ref{s1.1}. For a function $u: \bb R^d \to \bb R$ we define the {\em energy form} $\mc E(u,u)$ as
\[
\mc E(u,u) = \frac{1}{2} \iint  h(y-x) \big(u(y)-u(x)\big)^2 dx dy.
\]

For any two functions $u$, $v$ such that $\mc E(u,u) <+\infty$, $\mc E(v,v)<+\infty$, we define
\[
\mc E(u,v) = \frac{1}{2} \iint  h(y-x) \big(u(y)-u(x)\big)\big(v(y)-v(x)\big) dx dy.
\]

Notice that for functions $F$, $G$ in $\mc C_c^2(\bb R^d)$, $\mc E(F,G) = -\int F(x) \mc L G(x) dx$. 
Fix a reference density $\rho>0$. Define the {\em entropy} $\mc H:[0,\infty) \to [0,\infty)$ as
\begin{equation}
\label{ec1.1}
\mc H(a) = \int_\rho^a \log\Big(\frac{\phi(x)}{\phi(\rho)}\Big)dx.
\end{equation}

Fix $T>0$ and let $u_0: \bb R^d \to \bb R$ be a bounded function satisfying $\int |u_0(x)-\rho| dx <+\infty$. We say that a measurable function $u: \bb R^d \times [0,T] \to \bb R$ is an {\em energy solution} of the Cauchy problem
\begin{equation}
\label{ec2}
\left\{
\begin{array}{rcl}
\partial_t u &=& \mc L \phi(u) \\
u(0,\cdot) &=& u_0(\cdot)
\end{array}
\right.
\end{equation}
if:
\begin{itemize}
\item[i)] $u_t(\cdot)= u(\cdot,t)$ has {\em finite entropy}, that is,
\begin{equation}
\label{ec2.4}
\int_0^T \int \mc H(u(x,t))dxdt <+\infty,
\end{equation}
\item[ii)] $u_t$ satisfies the {\em energy estimate}
\begin{equation}
\label{ec2.3}
\int_0^T \mc E(\phi(u_t),\phi(u_t)) dt <+\infty,
\end{equation}
\item[iii)] for any smooth function $G: \bb R^d \times [0,T] \to \bb R$ with a compact support  contained in $\bb R^d \times [0,T)$,
\begin{equation}
\label{ec2.2}
\int_0^T \int \Big\{ \phi(u(x,t)) \mc L G(x,t)+  u(x,t)\partial_t G(x,t) \Big\}dx dt + \int G(x,0) u_0(x) dx =0.
\end{equation}
\end{itemize}

\begin{theorem}
\label{t2}
Let $u_0: \bb R^d \to [0,\infty)$ be a measurable, bounded initial profile with $\int (u_0(x)-u)^2dx <+\infty$ for some constant $u>0$ and let $\{\nu^n\}_n$ be a sequence of probability measures in $\Omega_{zr}$ associated to $u_0$. Let $\xi_t^n$ be the zero-range process $\xi_{t n^\alpha}$ starting from $\nu^n$ and let $\nu^n(t)$ be the distribution in $\Omega_{zr}$ of $\xi_t^n$. Assume that the interaction rate $g(\cdot)$ is non-decreasing and that there exist positive and finite constants $\rho$, $K$ such that
\begin{itemize}
\item[{\bf (H)}] For any $n \geq 0$, $H(\nu^n|\nu_\rho) \leq K n^d$.
\end{itemize}

If the interaction rate $g(\cdot)$ satisfies {\bf (FEM)}, also assume that there is a constant $\rho' >\rho$ such that 
\begin{itemize}
\item[{\bf (C)}] The measures $\nu^n$ are stochastically dominated by $\nu_{\rho'}$. 
\end{itemize}

Then $\{\nu^n(t)\}_n$ is associated to the function $u(\cdot,t)$, where $u(x,t)$ is the unique energy solution of (\ref{ec2}).
\end{theorem}

We say that equation (\ref{ec2}) is the {\em hydrodynamic limit} of the zero-range process $\xi_t$. 
Remember that the function $\phi(\rho)$ appearing in the hydrodynamic equation is equal to the expectation of the interaction rate $g(\cdot)$ with respect to the invariant measure $\nu_\rho$. For an integrable initial profile $u_0$, bounded in the complement of some ball around the origin, it is not difficult to see that the product measures $\nu^n$ defined in $\Omega_{zr}$ by substituting the Bernoulli marginals in (\ref{ec1.5}) by the measures $q_{u_0^n(z)}$ are associated to $u_0(\cdot)$. 
These measures also satisfy hypothesis {\bf (H)} if $\int |u_0(x)-\rho|dx <+\infty$. When {\bf (H)} is satisfied, we necessarily have $u=\rho$. Hypothesis {\bf (C)} is quite restrictive: it is satisfied by the measures $\nu^n$ defined above if and only if $\sup_x u_0(x) \leq \rho'$. In Appendix \ref{B} we explain how to get rid of  hypothesis {\bf (C)}. Perhaps the most restrictive assumption is $g(\cdot)$ being non-decreasing, since this is an assumption on the dynamics and not on the initial profile. We will see that the only point where we need $g(\cdot)$ to be non-decreasing is to obtain enough conditions on the limiting profile $u(x,t)$ in order to guarantee uniqueness of the hydrodynamic equation \eqref{ec2}. In Appendix \ref{B} we state a result which does not require $g(\cdot)$ to be non-decreasing, conditioned on a stronger uniqueness result for \eqref{ec2} which we conjecture to be true.

From now on, given a probability measure $\nu$ in $\Omega_{zr}$ ($\Omega_{ex}$ resp.) we denote by $\bb P_{\nu}$ the distribution of the process $\xi_\cdot^n$ ($\eta_\cdot^n$ resp.) starting from the initial distribution $\nu$ and we denote by $\bb E_{\nu}$ the expectation with respect to $\bb P_\nu$. For the invariant measures $\nu_\rho$, we will write $\bb P^\rho=\bb P_{\nu_\rho}$ and $\bb E^\rho = \bb E_{\nu_\rho}$.

\subsection{Entropy estimates}
In this Section we discuss the relevance of hypothesis {\bf (H)}. The main point is that the entropy $H(\nu^n(t)|\nu_\rho)$ is decreasing in time, and therefore it can be used as a Lyapunov function for the evolution of $\xi_t$. The results in this section are standard; they were introduced in \cite{GPV} and we include them here for the sake of completeness. We will follow the exposition of \cite{KL}.

First, we recall a variational formula for the relative entropy $H(\mu|\nu)$:
\[
H(\mu|\nu) = \sup_{f} \Big\{ \int f d\mu - \log \int e^f d\nu\Big\},
\]
where the supremum is over all  functions $f$ which are integrable with repect to $\mu$. A very useful way to estimate the integral of a function $f$ with respect to $\mu$ in terms of the relative entropy $H(\mu|\nu)$ is obtained taking $\gamma f$, $\gamma >0$ as a test function in the formula above:
\begin{equation}
\label{ec2.1}
\int f d\mu \leq \frac{H(\mu|\nu)}{\gamma} + \frac{1}{\gamma} \log \int e^{\gamma f} d\nu.
\end{equation}

This inequality is known as the {\em entropy inequality}. Of course, this inequality is not useful unless we have a good way to estimate relative entropies.
For $t \geq 0$, define $f_t^n = d\nu^n(t)/d\nu_\rho$. The density $f_t^n$ satisfies the Kolmogorov equation
\[
\frac{d}{dt} f_t^n = n^\alpha L_{zr}^* f_t^n.
\]

Notice that in our case the measure $\nu_\rho$ is reversible and therefore $L_{zr}^*=L_{zr}$, although this point is not crucial. Define $H_n(t) = H(\nu^n(t)|\nu_\rho)$. We see that
\begin{align*}
\frac{d}{dt} H_n(t) 
	& =\frac{d}{dt} \int f_t^n \log f_t^n d \nu_\rho = \int (1+\log f_t^n) n^\alpha L_{zr}^* f_t^n d\nu_\rho\\
	&= \int f_t^n n^\alpha L_{zr} (1+\log f_t^n) d\nu_\rho = n^\alpha \int f_t^n L_{zr} \log f_t^n d\nu_\rho. 
\end{align*}

The operator $L_{zr}$, being the generator of a particle system, is of the form $\sum_{i \in I} c_i(\xi) [f(\xi^i)-f(\xi)]$ for some set of indices $i \in I$, some non-negative rates $c_i$ and some transformations $\xi \mapsto \xi^i$ of the space $\Omega_{zr}$. Using the elementary inequality $a(\log b-\log a) \leq 2\sqrt a(\sqrt b - \sqrt a)$, on each one of the terms composing $f_t^n L_{zr} \log f_t^n$, we obtain 
\[
\frac{d}{dt} H_n(t) \leq -n^\alpha \sum_{x,y \in \bb Z^d} p(y-x) \int g(\xi(x)) \big[\sqrt{f_t^n(\xi^{x,y})} - \sqrt{f_t^n(\xi)} \big]^2 d\nu_\rho.
\]

Define $\mc D(f) = -\int \sqrt{f} L_{zr} \sqrt{f} d\nu_\rho$. Then,
\[
\frac{d}{dt} H_n(t) \leq -2 n^\alpha \mc D(f_t^n).
\]

Integrating the previous inequality between $0$ and $t$, we obtain that
\begin{equation*}
H_n(t) + 2n^\alpha \int_0^t \mc D(f_s^n) ds \leq H_n(0) \leq Kn^d,
\end{equation*}
where the last inequality is due to hypothesis {\bf (H)}. In other words, entropy is decreasing in time, and moreover it also controls the growth of the so-called {\em Dirichlet form} $\mc D(f_t^n)$. Define $\bar f_t^n  =t^{-1} \int_0^t f_s^n ds$. By convexity of $\mc D$, we conclude that
\begin{equation}
\label{ec3}
\mc D(\bar f_t^n) \leq \frac{K n^{d-\alpha}}{2t}.
\end{equation}

We can think about $H(\nu^n(t)|\nu_\rho)/n^d$ as a measure of the {\em macroscopic} entropy of the system. What is remarkable is that the bound (\ref{ec3}) on the Dirichlet form is enough to control the space-time fluctuations of the density of particles, as we will see in the following sections. This observation was introduced in \cite{GPV} and it is at the heart of the proof of hydrodynamic limits for particle systems.

\section{The exclusion process with long jumps}
\label{s2}
In this section we prove Theorem \ref{t1}. Let us define the {\em empirical process} $\pi_t^n$ by
\[
\pi_t^n(dx) = n^{-d} \sum_{z \in  \bb Z^d} \eta_t^n(z) \delta_{z/n}(dx),
\]
which turns out to be a Markov process in $\mc D([0,\infty),\mc M_+(\bb R^d))$. Notice that the weak topology in $\mc M_+(\bb R^d)$ is metrizable. Theorem \ref{t1} is an immediate consequence of the following result:

\begin{theorem}
\label{t3}
Under the hypothesis of Theorem \ref{t1}, the process $\{\pi_t^n; t \in [0,T]\}$ converges in distribution to the deterministic trajectory $u(x,t)dx$, where $u(x,t)$ is the solution of (\ref{ec1}).
\end{theorem}

A standard continuation argument shows that this theorem also holds for the process $\pi^n_\cdot$ defined in $\mc D([0,\infty),\mc M_+(\bb R^d))$. We restrict ourselves to a bounded interval to simplify some of the arguments.
The proof of this theorem follows the usual approach to convergence in distribution of stochastic processes. First we prove tightness of the distributions of $\{\pi^n_\cdot\}_n$. Then we prove uniqueness of the possible limiting points. Since a relatively compact sequence on a metrizable space with only one accumulation point is necessarily convergent, Theorem \ref{t3} follows from these two affirmations. 

\subsection{Tightness}
\label{s2.1}
We recall that a possible choice for a metric that generates the weak topology in $\mc M_+(\bb R^d)$ is
\[
d(\pi_1,\pi_2) =\sum_{i \geq 1} \frac{1}{2^i} \min\Big\{ \int G_id(\pi_1-\pi_2), 1\Big\}, 
\]
where $\{G_i\}_i$ is a numerable collection of suitable non-negative functions in $\mc C_c^\infty(\bb R^d)$, the set of infinitely differentiable functions with compact support. Denote $\pi_t^n(G) = \int G d\pi_t^n$. It is not difficult to see that the sequence $\{\pi_\cdot^n\}_n$ is tight if and only if the sequence $\{\pi_\cdot^n(G)\}$ is tight for any non-negative function $G \in \mc C_c^\infty(\bb R^d)$. Notice that the projections $\pi_t^n(G)$ are real-valued, and therefore easier to handle than the full process $\pi_\cdot^n$. By Dynkin's formula, 
\begin{equation}
\label{ec4}
\pi_t^n(G) = \pi_0^n(G) + \int_0^t \pi_s^n(\mc L_n G) ds + \mc M_t^n(G),
\end{equation}
where the operator $\mc L_n$ is defined by
\begin{equation}
\label{ec4.0}
\mc L_n  G(x/n) = \sum_{z \in \bb Z^d} n^\alpha p(z) \Big(G(x/n+z/n) - G(x/n)\Big)
\end{equation}
and $\mc M_t^n(G)$ is a martingale. The martingale $\mc M_t^n(G)$ has mean zero and quadratic variation
\[
\<\mc M_t^n(G)\> = \frac{1}{2}\int_0^t \frac{1}{n^{2d}} \sum_{y, z \in \bb Z^d} n^\alpha p(z-y) \Big(\eta_s^n(z)-\eta_s^n(y)\Big)^2\Big(G(z/n)-G(y/n)\Big)^2 ds.
\]

We will make repeated use of the identity $p(z)=h(z)=n^{-(d+\alpha)}h(z/n)$. Notice that
\[
\<\mc M_t^n(G)\> \leq \frac{t}{n^d} \sum_{y,z \in \bb Z^d} \frac{1}{2n^{2d}} h(z/n-y/n)\big(G(z/n)-G(y/n)\big)^2.
\]

This last sum is nothing but a Riemann sum for the  energy $\mc E(G,G)$. A Taylor expansion of $G$ shows that the integral defining $\mc E(G,G)$ is absolutely convergent, as well as the Riemann sum above. In particular, we have $\<\mc M_t^n(G)\> \leq C(G) t/n^d$. Therefore, the martingale $\mc M_t^n(G)$ converges to 0 in $\mc L^2(\bb P_{\mu^n})$. In the same way,
\[
\mc L_n G(x/n) = \frac{1}{n^d} \sum_{z \in \bb Z^d} h(z/n) \Big(G(z/n-x/n)-G(x/n)\Big)
\]
which is a Riemann sum for $\mc L G(x/n)$.
Since $G \in \mc C_c^\infty(\bb R^d)$ and due to the symmetry of $h$, it is not difficult to show that 
\begin{equation}
\label{ec4.1}
\lim_{n \to \infty} \sup_{z \in \bb Z^d} \big|\mc L_n G(z/n) - \mc LG(z/n)\big| =0,
\end{equation}
\begin{equation}
\label{ec4.2}
\lim_{n \to \infty} n^{-d} \sum_{z \in \bb Z^d} \big|\mc L_n G(z/n) - \mc LG(z/n)\big| =0.
\end{equation}

The simplest way to prove tightness of the sequence $\{\pi_\cdot^n(G)\}_n$ is to use Aldous' criterion, which now we explain.

\begin{proposition}[Aldous' criterion]
Let $(\mc E,d)$ be a separable, complete metric space.
Let $\{\bb P^n\}$ be a sequence of probability measures in $\mc D([0,\infty), \mc E)$. The sequence $\{\bb P^n\}_n$ is tight if:
\begin{itemize}
\item[i)] for any $\epsilon >0$ there exists a compact set $K \subseteq \mc E$ such that 
\[
\sup_n \bb P^n(\pi_0 \notin K) \leq \epsilon,
\]
\item[ii)] for any $\epsilon >0$ and any $T>0$,
\[
\lim_{\delta \to 0} \limsup_{n \to \infty} \sup_{\substack{\tau \in \mc T_T\\ \gamma \leq \delta}} 
	\bb P^n\big(d(\pi_{(\tau+\gamma) \wedge T}, \pi_\tau)>\epsilon \big) =0,
\]
where $\mc T_T$ is the set of stopping times bounded by $T$.
\end{itemize}
\end{proposition}

In our case, condition i) is automatically satisfied due to the fact that $\{\mu^n\}_n$ is associated to the bounded profile $u_0$. Condition ii) follows from equation (\ref{ec4}). In fact, since the number of particles per site is bounded by 1,
\begin{equation}
\label{ec5}
\Big| \int_{\tau}^{(\tau+\gamma)\wedge T} \pi_s^n(\mc L_n G) ds\Big| \leq \frac{\gamma}{n^d} \sum_{z \in \bb Z^d} \big|\mc L_n G(z/n)\big|.
\end{equation}

It is easy to see that $| \mc L_n G(x)| \leq C(G)/(1+||x||^{d+\alpha})$ for a constant $C(G)$ depending only on $||G||_\infty$, $||G''||_\infty$ and the support of $G$ (here $G''$ denotes the Hessian of $G$). Plugging this bound into the inequality (\ref{ec5}), condition ii) follows for the integral part. By the optional stopping theorem and Tchebyshev's inequality we have
\begin{align*}
\bb P_{\mu^n} \big(\big|\mc M_{(\tau+\gamma)\wedge T}^n(G) -\mc M_\tau^n(G)\big| >\epsilon\big)
	& \leq \epsilon^{-2} \bb E_{\mu^n} \big[\big(\mc M_{(\tau+\gamma)\wedge T}^n(G) - \mc M_\tau^n(G)\big)^2\big] \\
	&= \epsilon^{-2} \bb E_{\mu_n} \big[\<\mc M_{(\tau+\gamma)\wedge T}^n(G)\> - \<\mc M_\tau^n(G)\>\big]\\
	&\leq \frac{C(G) \gamma}{\epsilon^2 n^d},\\ 
\end{align*}
which goes to $0$ as $n \to \infty$, uniformly in $\tau$ and $\gamma \leq \delta$. In the last line above, the constant $C(G)$ depends only on $\|\nabla G\|_\infty$ and the support of $G$. Here and in the sequel we denote by $C$ a generic constant which may change from line to line, but depends only on the parameters indicated (for example, in the lines above $C(G)$ depends only on $G$). Therefore, the three terms on the right-hand side of (\ref{ec4}) are tight, from where tightness for $\{\pi_\cdot^n(G)\}_n$ (and in consequence for $\{\pi_\cdot^n\}_n$) follows.

\subsection{Uniqueness of limit points}
Once we have proved tightness for $\{\pi_\cdot^n\}$, we know that this sequence has accumulation points with respect to the topology of convergence in distribution. Let $\pi_\cdot$ be one of these points. Denote by $n'$ a subsequence for which $\pi_\cdot^{n'}$ converges to $\pi_\cdot$ The idea is to to take the limit through the subsequence $n'$ in (\ref{ec4}). By definition, $\pi_t^{n'}(G)$ converges to $\pi_t(G)$. By the assumptions on initial distributions, $\pi_0^{n'}(G)$ converges to $\int G(x) u_0(x) dx$. The martingale term $\mc M_t^n(G)$ converges to $0$ in $\mc L^2(\bb P_{\mu_n})$ and in particular it converges to $0$ also in distribution. However, $\mc L_n G$ is not a function in $\mc C_c(\bb R^d)$, due to the non-local character of the operator $\mc L_n$, so some justification is needed before taking the limit through $n'$ of $\pi_s^{n}(\mc L_n G)$. Observe that the number of particles per site is bounded by $1$. Using (\ref{ec4.1}), we can substitute $\mc L_n G$ by $\mc L G$ in (\ref{ec4}) by introducing an error term that vanishes as $n \to \infty$. Notice that
\[
\lim_{M \to \infty} n^{-d} \sum_{|z| \geq Mn} \big| \mc L G(z/n)\big| =0,
\]
uniformly in $n$. In particular, we can approximate $\mc LG$ by functions of compact support to obtain that
\[
\lim_{n' \to \infty} \int_0^t \pi_s^{n'} (\mc L G) ds = \int_0^t \pi_s(\mc L g) ds,
\]
in distribution. Therefore, the measure $\pi_\cdot$ satisfies
\[
\pi_t(G) = \pi_0(G) + \int_0^t \pi_s(\mc L G) ds
\]
for any function $G \in \mc C_c^\infty(\bb R^d)$. Let $G_t: [0,T] \to \mc C_c^\infty(\bb R^d)$ be a differentiable trajectory. Repeating the arguments above, it is not difficult to prove that $\pi_\cdot$ satisfies
\begin{equation}
\label{ec6}
\pi_t(G_t) = \pi_0(G_0) + \int_0^t \pi_s\big((\partial_t+\mc L) G\big) ds.
\end{equation}

Since the number of particles per site is bounded by $1$, it is not hard to see that the limiting measure $\pi_t$ is, for any $t >0$, absolutely continuous with respect to Lebesgue measure. Moreover, its density is bounded between $0$ and $1$. Let us write $\pi_t(dx) = u(x,t) dx$. In terms of the random (at this point) density $u(x,t)$, equation (\ref{ec6}) states that
\[
\int u(x,t) G_t(x) dx = \int u_0(x) G_0(x) dx + \int_0^t \int u(x,s)\big(\partial_t+ \mc L\big) G(x) dx ds
\]
for any smooth trajectory $G_t$, which is exactly the weak formulation of the hydrodynamic equation (\ref{ec1}). In other words, we have proved that $\pi_\cdot(dx) = u(x,\cdot)dx$ is concentrated on weak solutions of (\ref{ec1}). But this equation has only one weak solution starting from $u_0$ (see Section \ref{s7.1}). This uniqueness result finishes the proof of Theorem \ref{t3}. 

\section{The zero-range process with long jumps} 
\label{s3}
As we did for the exclusion process, we will consider the empirical measure 
\[
\pi_t^n(dx)= \frac{1}{n^d} \sum_{z \in \bb Z^d} \xi_t^n(z) \delta_{z/n}(dx)
\]
and we will prove that

\begin{theorem}
\label{t4}
Under the hypothesis of Theorem \ref{t2},
the sequence $\{\pi_\cdot^n\}_n$ is relatively compact with respect to the topology of convergence in distribution in the Skorohod space  $\mc D([0,T],\mc M_+(\bb R^d))$. All the limit points are concentrated on finite entropy solutions of the hydrodynamic equation (\ref{ec2}). When (\ref{ec2}) has a unique solution $u(x,t)$, $\pi_\cdot^n$ converges in probability to the deterministic path $u(x,\cdot)dx$.
\end{theorem}

We will explain what we mean by a finite entropy solution in Section \ref{s6}, where we state a more precise version of this Theorem.
The proof of this Theorem follows the same strategy followed in order to prove Theorem \ref{t3}, but it is technically more involved. For a proof in the case of a finite range, mean zero transition rate $p(\cdot)$, we refer to Chapter 5 of \cite{KL}. For the reader's convenience, we follow closely the proof in \cite{KL}, modifying to our setting.

\subsection{Some elementary estimates}
\label{s3.1}
Before we enter into the proof of Theorem \ref{t3}, in this section we collect some elementary estimates that will be used repeatedly. The results are well known, and we collect them here for the reader's convenience. Let $G \in \mc C_c^\infty(\bb R^d)$ be given. Define $F_0(x) = 1/(1+\|x\|^{d+\alpha})$. Our first estimate gives the behavior of $\mc L_n G$ and $\mc L G$ at infinity.

\begin{lemma}
\label{aux1}
For any function $G \in \mc C_c^\infty(\bb R^d)$, there exists a constant $C(G)$ such that for any $z \in \bb Z^d$ and any $x \in \bb R^d$,
\[
\big|\mc L_n G(z/n) \big| \leq C(G)F_0(z/n),
\]
\[
\big|\mc LG(x) \big| \leq C(G)F_0(x).
\]
\end{lemma}

\begin{proof}
Observe that, since $h(\cdot)$ is strictly positive, there exists a constant $\epsilon_0$ such that $\epsilon_0 \|x\|^{-(d+\alpha)} \leq h(x) \leq \epsilon_0^{-1} \|x\|^{-(d+\alpha)}$ for any $x \neq 0$. Therefore, it is enough to consider the case $h(x) = \|x\|^{-(d+\alpha)}$. For this last case, the result of the lemma follows easily from the compactness of the support of $G$ and a Taylor expansion of second order.
\end{proof}

The second estimate concerns the behavior of the moment generating function.

\begin{lemma}
\label{aux2} The function $\log M_\rho(\cdot)$ is strictly convex and increasing. In particular, if $M_\rho(\theta_0)$ is finite,
\[
\log M_\rho(\theta) \leq \frac{M_\rho(\theta_0) \theta}{\theta_0}
\]
for any $\theta \leq \theta_0$.
\end{lemma}

\begin{proof}
Let us denote by $q_{\rho,\theta}$ the distribution in $\bb N_0$ with density $e^{\theta n}/M_\rho(\theta)$ with respect to $q_\rho$. Let us denote by $E_\theta$ the expectation with respect to $q_{\rho,\theta}$. It is enough to observe that for $\Psi(\theta)=\log M_\rho(\theta)$, $\Psi'(\theta)=E_\theta[n]$ and $\Psi''(\theta) = \text{Var}_\theta[n]$.
\end{proof}

For a given function $G: \bb R^d \to \bb R$, let us define
\[
\|G\|_{1,n} = \frac{1}{n^d} \sum_{x \in \bb Z^d} \big| G(x/n)\big|.
\] 

\begin{lemma}
\label{aux3}
Let $G: \bb R^d \to \bb R$ be a bounded function. There exists a constant $C>0$ such that for any $t \geq 0$ we have
\[
\bb E_{\nu^n} \big[ \big| \pi_t^n(G)\big|\big] \leq C\big( \|G\|_\infty + \|G\|_{1,n}\big).
\]
\end{lemma}

\begin{proof}
Without loss of generality, we can assume that $G$ is non-negative. Remember that $\bb E_{\nu^n}[\pi_t^n(G)] = \int \pi^n(G) f_t^n d\nu_\rho$. By the entropy inequality,
\begin{align*}
\bb E_{\nu^n}[\pi_t^n(G)]
	&\leq \frac{K}{\gamma} + \frac{1}{\gamma n^d} \log \int \exp\Big\{ \gamma\sum_{z \in \bb Z^d} \xi(z) G(z/n)\Big\}d \nu_\rho\\
	& \leq \frac{K}{\gamma} + \frac{1}{\gamma n^d} \sum_{z \in \bb Z^d} \log M_\rho(\gamma G(z/n)).
\end{align*}
Take $\theta_0>0$ such that $M_\rho(\theta_0)$ is finite. Then take $\gamma = \theta_0/\|G\|_\infty$. Since the measure $\nu_\rho$ is of product form, using Lemma \ref{aux2} we can bound $\bb E_{\nu^n}[\pi_t^n(G)]$ by
\[
\frac{K\|G\|_\infty}{\theta_0}+ \frac{\|G\|_\infty}{\theta_0 n^d} \sum_{z \in \bb Z^d} \frac{M_\rho(\theta_0)G(z/n)}{\|G\|_\infty}.
\]
Taking $C= \theta_0^{-1}\max\{K, M_\rho(\theta_0)\}$ we obtain the desired bound.
\end{proof}

\subsection{Tightness}
Remember that in order to prove tightness for $\{\pi_\cdot^n\}_n$, it is enough to prove tightness for $\{\pi_\cdot^n(G) \}_n$ for any non-negative $G \in \mc C_c^\infty(\bb R^d)$. The projection of the empirical measure can be written as
\begin{equation}
\label{ec7}
\pi_t^n(G) = \pi_0^n(G) + \int_0^t \frac{1}{n^d} \sum_{z \in \bb Z^d} g\big(\xi_s^n(x)\big) \mc L_n G(z/n) ds + \mc M_t^n(G),
\end{equation}
where $\mc M_t^n(G)$ is a martingale of quadratic variation 
\[
\<\mc M_t^n(G) \> = \int_0^t \frac{1}{n^{3d}} \sum_{y \in \bb Z^d} g(\xi_s^n(y)) \sum_{z \in \bb Z^d} h(z/n-y/n) \big(G(z/n)-G(y/n)\big)^2 ds.
\]

Let us define
\[
\mc Q_n G(x/n) = \frac{1}{n^d} \sum_{z \in \bb Z^d} h(z/n-x/n) \big(G(z/n)-G(x/n)\big)^2.
\]

Notice that $\mc Q_n G(x/n)$ is  a Riemann sum for the integral $\int h(y-x) (G(y)-G(x))^2dy$ evaluated at $x/n$. We can rewrite the quadratic variation of $\mc M_t^n(G)$ as 
\[
\<\mc M_t^n(G) \> = \int_0^t \frac{1}{n^{2d}} \sum_{z \in \bb Z^d} g(\xi_s^n(z)) \mc Q_n G(z/n) ds.
\]

In this expression, as for the exclusion process with long jumps, there is an extra $1/n^d$ term that will make the quadratic variation of $\mc M_t^n(G)$ converge to 0. Now the difference is that we need to bound the expectation of $g(\xi_t^n(x))$, which is not longer a bounded random variable. We can write $\<\mc M_t^n(G)\>$ in terms of $\bar f_t^n$:
\[
\<\mc M_t^n(G) \> =  \frac{t}{n^d} \int \frac{1}{n^d} \sum_{z \in \bb Z^d} g(\xi(z)) \mc Q_n G(z/n) \bar f_t^n d\nu_\rho.
\]

Remember the bound $g(n) \leq \kappa n$. Notice that the entropy of $\bar f_t^n$ is also bounded by $Kn^d$. Therefore, the proof of Lemma \ref{aux3} also applies here, so we obtain the bound
\[
\bb E_{\nu^n}\<\mc M_t^n(G) \> \leq \frac{Ct}{n^d}\big(\|\mc Q_nG\|_\infty + \|\mc Q_nG\|_{1,n}\big).
\]
Notice that both norms $\|\mc Q_nG\|_\infty + \|\mc Q_nG\|_{1,n}$ are uniformly bounded in $n$.
We conclude that $\mc M_t^n(G)$ converges to $0$ in $\mc L^2(\bb P_{\nu^n})$ as $n$ goes to $\infty$, like in the case of the exclusion process. Therefore, $\{\mc M_\cdot^n(G)\}$ is tight.

With the same notation of Section \ref{s2.1} and using the entropy estimate, for any non-negative, bounded function $F$ we have
\begin{align*}
\bb E_{\nu^n} \Big[ \int_\tau^{(\tau+\gamma)\wedge T} 
	&\frac{1}{n^d} \sum_{z \in \bb Z^d} g(\xi_s^n(z)) F(z/n) ds \Big] \leq \frac{K}{\gamma_0} + \\
	&+\frac{1}{\gamma_0 n^d} \log \bb E^\rho \Big[\exp\Big\{\gamma_0 \int_\tau^{(\tau+\gamma)\wedge T} \sum_{z \in \bb Z^d} g(\xi^n_s(z)) F(z/n)ds\Big\}\Big].\\
\end{align*}

By Jensen's inequality, this last line is bounded by
\begin{multline*}
\leq  \frac{K}{\gamma_0} + \frac{1}{\gamma_0 n^d} \log \bb E^\rho \Big[\frac{1}{\gamma} \int_\tau\limits^{(\tau+\gamma)\wedge T} \exp\Big\{ \gamma_0 \gamma \sum_{z \in \bb Z^d} g(\xi^n_s(z))) F(z/n) \Big\} ds\Big]\\
	\leq \frac{K}{\gamma_0} + \frac{1}{\gamma_0 n^d} \log \bb E^\rho\Big[ \frac{1}{\gamma} \int_0^T \exp\Big\{ \gamma_0 \gamma \sum_{z \in \bb Z^d} \kappa \xi^n_s(z) F(z/n)\Big\}ds\Big],
\end{multline*}
where we have used the fact that $g(n) \leq \kappa n$ in the last line. Since the measure $\nu_\rho$ is invariant under the evolution of $\xi_t^n$, the expectation is bounded by
\begin{multline*}
\frac{K}{\gamma_0} + \frac{1}{\gamma_0 n^d} \log \Big\{ \frac{T}{\gamma} \prod_{z \in \bb Z^d} M_\rho(\gamma_0 \gamma \kappa F(z/n) )\Big\}=\\
 = \frac{K}{\gamma_0} + \frac{1}{\gamma_0 n^d}\Big\{ \log(T/\gamma) + \sum_{z \in \bb Z^d} \log M_\rho(\gamma_0 \gamma \kappa F(z/n))\Big\}. 
\end{multline*}
For $\gamma$ small enough (how small depends only on $\|F\|_\infty$), by Lemma \ref{aux2} the expression above in bounded by 
\[
\frac{K}{\gamma_0} + \frac{1}{\gamma_0 n^d}\Big\{ \log(T/\gamma) + \sum_{z \in \bb Z^d} C(\rho) (\theta_0) \gamma_0 \gamma \kappa F(z/n) \Big\},
\]
where the last inequality is true for $\gamma$ small enough and $C(\rho)$ si simply $M_\rho(\theta_0)/\theta_0$. We conclude that
\[
\sup_{\substack{\gamma \leq \delta\\ \tau \in \mc T_T}} \limsup_{n \to \infty} \bb P_{\mu^n} \Big( \int\limits_\tau^{(\tau+\gamma)\wedge T} 
	\frac{1}{n^d} \sum_{z \in \bb Z^d} g(\xi_s^n(z)) F(z/n) ds>\epsilon\Big)
		\leq  \frac{K}{\gamma_0} + C(\rho)\delta \kappa \int F(x) dx
\]
for $\delta$ small enough, if $F$ is such that the Riemann sum above converges as $n \to \infty$. For $G \in \mc C_c^\infty(\bb R^d)$ and using Lemma \ref{aux1}, tightness follows for the integral term in (\ref{ec7}) by Aldous' criterion. This finishes the proof of tightness for the sequence $\{\pi_\cdot^n(G)\}_n$. 

\subsection{Identification of limit points}
\label{s3.2}
Let $G: \bb R^d \times [0,\infty) \to \bb R$ be a function of class $\mc C^{2,1}$ and bounded support. Write $G_t(\cdot) = G(\cdot,t)$. Following the same computations leading to tightness made on the previous section, we see that
\begin{multline}
\label{ec8}
\pi_t^n(G_t) = \pi_0^n(G_0) + \int_0^t \pi_s^n(\partial_t G_s) ds +\\
	+ \int_0^t \frac{1}{n^d} \sum_{z \in \bb Z^d} g(\xi_s^n(z)) \mc L_n G_s(z/n) ds + \mc M_t^n(G), 
\end{multline}
where $\mc M_t^n(G)$ is a martingale such that $\bb E_{\nu^n}[\mc M_t^n(G)^2] \to 0$ as $n \to \infty$. We already know that the sequence $\{\pi_\cdot^n\}_n$ is tight, so now we want to characterize its limit points as energy solutions of the hydrodynamic equation (\ref{ec2}). The main difference between (\ref{ec4}) and (\ref{ec8}) is the presence of the function $g(\xi_s^n(x))$. This makes the second integral in (\ref{ec8}) not to be a function of the empirical measure $\pi_\cdot^n$. Our objective, then, is to write the integral 
\[
 \int_0^t \frac{1}{n^d} \sum_{z \in \bb Z^d} g(\xi_s^n(z)) \mc L_n G_s(z/n) ds
\]
as a function of the empirical measure plus a rest that vanishes in distribution as $n \to \infty$. This is the content of the following theorem, known in the literature as the Replacement Lemma:

\begin{theorem}[Replacement Lemma]
\label{t5}
Under the hypothesis of Theorem \ref{t2}, for any $M >0$ and any $t >0$,
\begin{multline*}
\lim_{\epsilon \to \infty} \limsup_{n \to \infty} \bb E_{\nu^n} \Big[ \int_0^t \frac{1}{n^d} \sum_{|x| \leq Mn}\Big| \frac{1}{(2\epsilon n +1)^d} \sum_{|y| \leq \epsilon n} g(\xi_s^n(x+y))\\
	 -\phi\Big(\frac{1}{(2\epsilon n +1)^d} \sum_{|y| \leq \epsilon n} \xi_s^n(x)\Big) \Big| ds\Big] =0.
\end{multline*}
\end{theorem}

Here and below, by abuse of notation we write $\epsilon n$ in place of its integer part $[\epsilon n] = \sup\{n \in \bb Z; n\leq \epsilon n\}$. Notice that with this convention, $(2\epsilon n +1)^d$ is just the cardinality of the set $\{|y| \leq \epsilon n\}$. We postpone the proof of this Theorem to the next section. The Replacement Lemma is, with no doubt, the heart of the proof of Theorem \ref{t4}. Let $\pi_\cdot$ be a limit point of $\{\pi_\cdot^n\}_n$, and denote by $n'$ a subsequence for which $\pi_\cdot^{n'} \to \pi_\cdot$ in distribution. Notice that
\[
\frac{1}{(2 \epsilon n+1)^d} \sum_{|z|\leq \epsilon n} \xi_s^n(x+z) = c_n \pi_s^n\big((2\epsilon)^{-d}
	\mathbf{1}(|y-x/n| \leq \epsilon)\big),
\]
where $c_n$ is a normalizing constant that goes to 1 as $n \to \infty$. In particular, this expression is a function of the empirical measure. Define, for any measure $\pi \in \mc M_+(\bb R^d)$, the function $I_\epsilon \pi$ by
\[
I_\epsilon \pi(x) = (2\epsilon)^{-d} \int_{|y-x| \leq \epsilon} \pi(dy).
\]

Since the function $\mc L_n G_s(x)$ is uniformly continuous and bounded by $C(G) F_0(x)$, thanks to the Replacement Lemma we can write
\begin{multline*}
\int_0^t \frac{1}{n^d} \sum_{z \in \bb Z^d} g(\xi_s^n(z)) \mc L_n G(x/n,s) ds =\\
	= \int_0^t \frac{1}{n^d} \sum_{z \in \bb Z^d} \phi\big(I_\epsilon \pi_s^n(x/n)\big) \mc L_n G(x/n,s) ds + R_t^{n,\epsilon} (G),
\end{multline*}
where $R_t^{n,\epsilon} (G)$
is an error term that vanishes in $\mc L^1(\bb P_{\nu^n})$ when $n \to \infty$ and then $\epsilon \to 0$. We denote by $R_t^{n,\epsilon}(G)$ any term with this property. 
	
Using (\ref{ec4.1}) and (\ref{ec4.2}), we can write
\begin{multline*}
\int_0^t \frac{1}{n^d} \sum_{z \in \bb Z^d} \phi\big(I_\epsilon \pi_s^n(x/n)\big) \mc L_n G(x/n,s) ds=\\
	= \int_0^t \int _{\bb R^d} \phi\big(I_\epsilon \pi_s^n(x)\big) \mc L G(x,s) dx ds + R_t^{\epsilon,n}(G).
\end{multline*}

Notice that we can not say that $I_\epsilon \pi_s^{n'}(x) \to I_\epsilon \pi_s(x)$, since the indicator function $\mathbf{1}(|\cdot-x| \leq \epsilon)$ is not a continuous function. However, by Portmanteau's lemma, this is true whenever $\pi_s(\partial\{y; |y-x| \leq \epsilon\}) =0$, where we denote by $\partial A$ the boundary of $A \subseteq \bb R^d$.  Since $\pi$ is a Radon measure, this is the case in a set of full measure in $\bb R^d$. In particular, $I_\epsilon \pi_s^{n'}(x) \to I_\epsilon \pi_s(x)$ for $x$ in a set of full measure in $\bb R^d$. Since $| \phi\big(I_\epsilon \pi_s^n(x)\big) \mc L G(x,s)|$ is bounded by $ \kappa I_\epsilon \pi_s^n(x) |\mc L G(x,s)|$, by the dominated convergence theorem we conclude that
\[
\lim_{n' \to \infty} \int_0^t \int_{\bb R^d} \phi\big( I_\epsilon \pi_s^{n'} (x) \big) \mc L G(x,s) dx ds 
	= \int_0^t \int_{\bb R^d} \phi\big( I_\epsilon \pi_s(x)\big) \mc L G(x,s) dx ds.
\]

Therefore, taking the limit in (\ref{ec8}) through the subsequence $n'$, we obtain that
\begin{align*}
\int G(x,t) \pi_t(dx) 
	&= \int u_0(x) G(x,0) dx + \int_0^t \int \partial_t G(x,s) \pi_s(dx) ds \\
	&+ \int_0^t \int \phi\big( I_\epsilon \pi_s(x)\big) \mc L G(x,s) dx ds + R_t^\epsilon(G),
\end{align*}
where $R_t^\epsilon(G)$ is a rest that vanishes in expectation when $\epsilon \to 0$. At this point, we only need to find the limit as $\epsilon \to 0$ of
\[
\int_0^t \int \phi\big( I_\epsilon \pi_s(x)\big) G(x,s) dx ds.
\]

It is clear that $I_\epsilon$ is just an approximation of the identity, so it is reasonable to expect that $I_\epsilon \pi_s \to \pi_s$ as $\epsilon \to 0$ in some sense. We need some regularity for $\pi_s$ in order to pass to the limit inside the function $\phi$. This is the case if, for example, $\pi_s$ is absolutely continuous with respect to Lebesgue measure. The following lemma says that this is indeed the case.

\begin{lemma}
\label{l1}
The process $\pi_\cdot$ is concentrated in trajectories of the form $u(x,\cdot)dx$, where $u:[0,T] \times \bb R^d \to \bb R$ is non-negative and locally integrable. Moreover,
\begin{equation}
\label{ec8.2}
\int_0^T \int \mc H(u(x,t)) dx dt <+\infty.
\end{equation}
\end{lemma}

We will asume this lemma and we will prove it later. Since $\pi_s(dx) = u(x,s) dx$ and $u(x,s)$ is locally integrable, we see that $I_\epsilon \pi_s(x)$ converges to $u(x,s)$ $a.s.$ with respect to Lebesgue measure,  and also in $\mc L_{\text{loc}}^1(\bb R^d)$. The bound $g(n) \leq Kn$ implies that $\phi(\rho) \leq K \rho$. This linear bound plus the fact that $\phi(\rho)$ is locally Lipschitz (since it is smooth) allow to conclude that $\phi(I_\epsilon\pi_s(x))dx$ converges weakly to $\phi(u(x,s))dx$ as $\epsilon \to 0$.
By Fubini-Tonelli's theorem, we conclude that
\[
\lim_{\epsilon \to 0} \int_0^t \int \phi\big( I_\epsilon \pi_s(x)\big) \mc LG(x,s) dx ds = \int_0^t \int \phi(u(x,s))dx ds.
\]

Therefore, $u(x,t)$ satisfies conditions i) and iii) of the definition of energy solutions of (\ref{ec2}). In order to finish the proof of Theorem \ref{t4} it is only left to prove that $u(x,t)$ satisfies the energy estimate $\int_0^T \mc E(\phi(u_t),\phi(u_t)) dt <+\infty$. This is the content of the following theorem:

\begin{theorem}[Energy estimate]
\label{t5.1}
Under the hypothesis of Theorem \ref{t2}, the limit points of $\{\pi_\cdot^n\}_n$ are concentrated in measures of the form $u(x,\cdot)dx$, where $u_t(\cdot)=u(\cdot,t)$ satisfies
\begin{equation}
\label{ec8.1}
\int_0^T \mc E(\phi(u_t),\phi(u_t)) dt <+\infty.
\end{equation}
\end{theorem}

Accepting the validity of Theorems \ref{t5}, \ref{t5.1} and Lemma \ref{l1}, we conclude that the limit point $\pi_\cdot$ is concentrated on enegy solutions of the hydrodynamic equation (\ref{ec2}). We will dedicate the next sections to the proof of Theorems \ref{t5}, \ref{t5.1}. We end this section proving Lemma \ref{l1}.

\begin{proof}[Proof of Lemma \ref{l1}]

The reader may recognize in the following lines an application of Varadhan's lemma. 
Let us  define $h(\theta) = \log M_\rho(\theta)$. For any sequence of random variables $a_1,\dots,a_k$ we have the estimate
\begin{equation}
\label{ec9}
\begin{split}
\log E[\exp\{\max_i\{a_i\}\}] 
	&= \log E[\max_i\{\exp\{a_i\}\}] \leq \log E\big[\sum_i \exp\{a_i\}\big]\\
	&= \log \sum_i E[\exp\{a_i\}] \leq \log \big\{k \max_i E[\exp\{a_i\}]\big\}\\
	&\leq \log k +\max_i \log E[\exp\{a_i\}]. 
\end{split}
\end{equation}

Let $G^1,\dots,G^k$ be a sequence of functions in $\mc C_c(\bb R^d \times [0,T])$. We will apply the estimate above to the sequence $\int_0^T\pi_t^n(G^1_t)dt,\dots,\int_0^T\pi_t^n(G_t^k)dt$.  For a function $G:\bb R^d \to \bb R$, let us define $J_n(G) = n^{-d}\sum_x h(G(x/n))$. 
By the entropy inequality,
\begin{align*}
\bb E_{\nu^n}\Big[ 
	&\quad
	\max_{i=1,\dots,k} 
	\Big\{ \frac{1}{T}\int_0^T\big(\pi_t^n(G_t^i) -  J_n(G^i_t)\big)dt\Big\}\Big] \leq \\
	&\quad
	\leq K + \frac{1}{n^d} \log \bb E^\rho\Big[\exp\Big\{\max_{i=1,\dots,k} \frac{n^d}{T} \int_0^T\big(\pi^n(G_t^i) -  J_n(G_t^i)\big)dt\Big\}\Big]\\
	&\quad
	\leq K +\frac{\log k}{n^d} +\max_{i=1,\dots,k} \frac{1}{n^d}\log \bb E^\rho\Big[ \exp\Big\{\frac{n^d}{T} \int_0^T\big(\pi^n(G_t^i) - J_n(G^i_t)\big)dt\Big\}\Big]\\
	&\quad
	\leq K +\frac{\log k}{n^d} +\max_{i=1,\dots,k} \frac{1}{n^d}\log \frac{1}{T}\int_0^T \bb E^\rho\Big[ \exp\big\{n^d \big(\pi^n(G_t^i) - J_n(G^i_t)\big)\big\}\Big]dt.\\
\end{align*}

Taking the limit along the subsequence $n'$, the left-hand side of this inequality converges to $E[\max_i T^{-1}\int_0^T \{\pi_t(G_t^i) - J(G_t^i)\}dt]$, and the right-hand side converges to $K$. Now $k$ is arbitrary, so we have proved that
\begin{equation}
\label{ec9.1}
E\Big[ \sup_G  \int_0^T \big\{\pi_t(G_t) - J(G_t)\big\}dt\Big] \leq K T,
\end{equation}
where the supremum is over ${G \in \mc C_{c,T}=:\mc C_c(\bb R^d\times [0,T])}$.
The supremum inside the integral can be computed. Denote by $\mc H$ the Legendre transform of $h$: $\mc H(a) = \sup_{\theta} \big\{ a \theta - h(\theta)\big\}$. 
By Lemma \ref{aux2}, $h(\cdot)$ is strictly convex and smooth. Therefore, $\mc H(\cdot)$ is also strictly convex, and grows at least linearly. Let us define
\[
\mc J(\pi) =\sup_{G \in \mc C_{c,T}}  \int_0^T \big\{\pi_t(G_t) - J(G_t)\big\}dt.
\]

Suppose that $\pi_t(dx)$ is absolutely continuous with respect to Lebesgue measure. Let $u(x,t)$ be its density. Then,
\begin{align*}
\mc J(\pi) 
	&= \sup_{G \in \mc C_{c,T}} \int_0^T\int \big\{ u(x,t) G(x) - h(G(x))\big\}dxdt\\
	&\leq \int_0^T \int \sup_{g \in \bb R} \big\{ u(x,t) g -h(g)\big\}dxdt= \int_0^T \int \mc H(u(x,t))dxdt.
\end{align*}

By Lemma \ref{aux2}, $h(\cdot)$ is strictly convex and smooth. Therefore, in the definition of $\mc H(a)$ the minimizer $\theta(a)$ such that $\mc H(a) = a\theta(a) - h(\theta(a))$ is unique and smooth as a function of $a$. In fact $\theta =\theta(a)$ satisfies $a=h'(\theta)$. After some computations, we can see that the function $\mc H(a)$ is equal to the entropy function defined in \eqref{ec1.1}.
We conclude that the supremum in the definition of $\mc J(\pi)$ is attained at the continuous function $\theta(u(x,t))$, and therefore $\mc J(\pi)= \int \mc H(u(x,t)) dx$. 
Now assume that $\pi_t$ is not absolutely continuous with respect to Lebesgue measure. We will prove that $\mc J(\pi)=+\infty$. 
Let $F$ be a bounded, closed set such that $\int_0^T \pi_t(F)dt=a>0$ and $\int_0^T\int_Fdxdt=0$. Fix $\beta>0$. Let $F_\epsilon=\{(y,t) \in \bb R^d \times [0,T];d((y,t),F)\leq \epsilon\}$ be the ball of radius $\epsilon$ and center $F$. 
Let $G^\epsilon: \bb R^d \times [0,T] \to [0,\beta]$ be a continuous function such that $G^\epsilon(x,t)=\beta$ if $(x,t) \in F_{\epsilon/2}$ and $G^\epsilon(x,t)=0$ if $(x,t) \notin F_\epsilon$. Then,
\[
\int_0^T \pi_t(G^\epsilon_t)dt  -\int_0^T J(G^\epsilon_t)dt  \geq a\beta - h(\beta)\int_{F_\epsilon}dxdt.
\]
Sending $\epsilon \to 0$, we conclude that $\mc J(\pi) \geq a \beta$. Since $\beta$ is arbitrary, we deduce that $\mc J(\pi)=+\infty$. By \eqref{ec9.1}, we know that $E[\mc J(\pi)]\leq K$. In particular, $\mc J(\pi)$ is finite {\em a.s.}, from where $\pi_\cdot$ is concentrated in trajectories of the form $u(x,\cdot)dx$. 
Since $h(\cdot)$ is strictly convex and smooth, $\mc H(\cdot)$ is also strictly convex, and grows at least linearly.
We conclude that  $u(x,t)$ is locally integrable.
\end{proof}

\section{The Replacement Lemma}
\label{s4}
In this section we prove Theorem \ref{t5}. Following \cite{GPV}, the proof is divided into two pieces: the so-called {\em one-block} and {\em two-blocks} estimates. The one-block estimate does not pose a real challenge for the model we are considering. After proving the one-block estimate, the two-blocks estimate follows from the {\em moving particle lemma}, which roughly states that the cost of moving a particle from one site $x$ to another site $y$ can be estimated by the Dirichlet form of the process restricted to a box containing both points. In the diffusive case, the moving particle lemma is an elementary application of Cauchy-Schwartz inequality, but in our superdiffusive setting, a more sophisticated proof is needed. For the sake of completeness, we give a proof of both one-block and two-blocks estimates. 

To simplify the notation, for $x \in \bb Z^d$ and $l>0$ we define
\[
\xi^l(x) = \frac{1}{(2l +1)^d} \sum_{|y| \leq l} \xi(x+y),
\]
\[
V_x^l(\xi) = \Big|  \frac{1}{(2l +1)^d} \sum_{|y| \leq l} g(\xi(x+y)) - \phi(\xi^l(x))\Big|.
\]

In terms of $V_x^l$, the Replacement Lemma can be written as
\begin{equation}
\label{ec12}
\lim_{\epsilon \to 0} \limsup_{n \to \infty} \bb E_{\nu^n} \Big[ \int_0^t \frac{1}{n^d} \sum_{|x|\leq Mn} V_x^{\epsilon n} (\xi_s^n) ds \Big] =0.
\end{equation}

The expectation in (\ref{ec12}) can be written as
\[
t\int \frac{1}{n^d} \sum_{|x|\leq Mn} V_x^{\epsilon n} (\xi) \bar f_t^n(\xi) d\nu_\rho. 
\]

It will be convenient to introduce an intermediate scale $l$ into (\ref{ec12}).
The limit in (\ref{ec12}) is a simple consequence of the following two lemmas (see \cite{KL}, Sect. 5.3):

\begin{lemma}[1-block]
\label{1b}
Under the hypothesis of Theorem \ref{t2},
\[
\lim_{l \to \infty} \limsup_{n \to \infty} \int \int \frac{1}{n^d} \sum_{|x|\leq Mn} V_x^{l} (\xi) \bar f_t^n(\xi) d\nu_\rho =0.
\]
\end{lemma}

\begin{lemma}[2-blocks]
\label{2b}
Under the hypothesis of Theorem \ref{t2},
\begin{equation}
\label{ec12.1}
\lim_{l \to \infty} \limsup_{\epsilon \to \infty} \limsup_{n \to \infty} \sup_{y,z} \int \frac{1}{n^d} \sum_{|x|\leq Mn} \big|\xi^l(x+z)-\xi^l(x+y)\big| \bar f_t^n(\xi) d\nu_\rho =0,
\end{equation}
where the supremum is over $y, z \in \bb Z^d$ such that $|y|, |z| \leq \epsilon n$, $|z-y| \geq 2l+1$.
\end{lemma}

Observe that ``pasting'' boxes of size $n$ in a convenient way, we can assume without loss of generality that $M=1$. Therefore, from now on we assume  $M= 1$. 

\subsection{The one-block estimate}

In this section we prove Lemma \ref{1b}. The proof is standard and we basically repeat the proof in Section 5.4 of \cite{KL}. We include a complete proof in order to simplify the exposition of the proof of Lemma \ref{2b}, for which new arguments are needed.

Remember that the density $\bar f_t^n$ satisfies the estimates $\int \bar f_t^n \log \bar f_t^n d\nu_\rho \leq K n^d$, $\mc D(\bar f_t^n) \leq K n^{d-\alpha}/2t$ and $\bar f_t^n \nu_\rho \preceq \nu_{\rho'}$. According to Section 5 of \cite{KL}, this is all we need to know about $\bar f_t^n$. In fact, the coupling estimate $\bar f_t^n \nu_\rho \preceq \nu_{\rho'}$ is only needed to prove the two-blocks estimate when $g(\cdot)$ satisfies {\bf (B)}.
First we introduce a cut-off of large densities. As in step 1, Section 5.4 of \cite{KL} when $g(\cdot)$ satisfies {\bf (B)} or Lemma 5.4.2 of \cite{KL} when $g(\cdot)$ satisfies {\bf (FEM)}, the following limit holds:
\[
\lim_{a \to \infty} \limsup_{l \to \infty} \limsup_{n \to \infty} 
	\int \frac{1}{n^d} \sum_{|x| \leq n} V_x^l(\xi) {\mathbf 1}\{\xi^l(x) \geq a\} \bar f_t^n d\nu_\rho =0.
\]
Let us define $V_x^{l,a}(\xi) = V_x^l(\xi) {\mathbf 1}\{\xi^l(x) <a\}$. In order to prove Lemma \ref{1b}, we just need to prove that for any $a>0$,
\begin{equation}
\label{ec13}
\limsup_{l \to \infty} \limsup_{n \to \infty} 
	\int \frac{1}{n^d} \sum_{|x| \leq n} V_x^{l,a}(\xi) \bar f_t^n d\nu_\rho =0
\end{equation}
From now on we simply write $\bar f_t^n = f$, keeping in mind the properties of $\bar f_t^n$ stated above. Let $\tau_x$ denote the translation by $x$: $\tau_x \xi(z) = \xi(x+z)$. Performing the change of variables $\bar \xi = \tau_x \xi$ and by translation invariance of the measure $\nu_\rho$, we can pass the sum in $x$ from $V_x^l$ to $f$. In this way we see that the integral in (\ref{ec13}) is equal to
\begin{equation}
\label{ec14}
\Big(\frac{2n+1}{n}\Big)^d \int V_0^{l,a} (\xi) \frac{1}{(2n+1)^d} \sum_{|x| \leq n} \tau_x f(\xi) d \nu_\rho,
\end{equation}
where we write $\tau_x f (\xi) = f(\tau_x \xi)$. The prefactor $((2n+1)/n)^d$ is bounded and it  can be ignored. Notice that the weighted sum in (\ref{ec14}) is also a density. We denote this density by $\bar f$. Define $\Lambda_l =\{-l,\dots,l\}^d$ and let $\mc F_l$ be the $\sigma$-algebra generated by $\{\xi(x); x \in \Lambda_l\}$. For simplicity, we also denote by $\Lambda_l$ the cardinality of $\{-l,\dots,l\}^d$, which is equal to $(2l+1)^d$. Notice that $\Lambda_n = \{x \in \bb Z^d; |x| \leq n\}$. Define $\bar f_l = E[\bar f|\mc F_l]$, that is, $\bar f_l$ is the conditional expectation of $\bar f$ with respect to the configuration $\xi$ restricted to the box $\Lambda_l$. With this notation, and due to the product structure of $\nu_\rho$, the integral in (\ref{ec14}) is equal to
\begin{equation}
\label{ec15}
\int V_0^{l,a}(\xi) \bar f_l (\xi) d\nu_\rho,
\end{equation}
where now the integral is over the set $\Omega_l = \bb N_0^{\Lambda_l}$. At this point we have reduced the original problem into a finite-dimensional, static problem. In fact, due to the indicator function, the integrand is different from zero only on the finite set $\{\xi \in \Omega_l ; \xi^l(0) <a\}$.

The following step is to estimate the Dirichlet form of $\bar f_l$. For any $x,y  \in \bb Z^d$, let us define
\[
\mc D^{x,y}(f) = p(y-x) \int g(\xi(x)) \Big(\sqrt{f(\xi^{x,y})} -\sqrt{f(\xi)}\Big)^2 d\nu_\rho.
\]

Notice that $\mc D^{x,y}(f)= \mc D^{y,x}(f)$, due to the symmetry of $p(\cdot)$.
The form $\mc D^{x,y}(f)$ is convex as a function of $f$. Since $\bar f$ is an average over translations of $f$, we have
\[
\mc D^{x,y}(\bar f) 
	=\frac{1}{\Lambda_n} \sum_{z \in \Lambda_n} \mc D^{x,y} (\tau_z f) = \frac{1}{\Lambda_n} \sum_{z \in \Lambda_n} \mc D^{x+z,y+z}(f).
\]

For a given density $f$, let us define
\[
\mc D_l(f) = \sum_{x,y \in \Lambda_l} \mc D^{x,y} (f).
\]

The form $\mc D_l(\cdot)$ corresponds to the Dirichlet form associated to the generator $L_{zr}^l $ of the process restricted to the box $\Lambda_l$. We have
\[
\mc D_l(\bar f) \leq \sum_{\substack{x,y \in \Lambda_l\\ z \in \Lambda_n}} \frac{1}{\Lambda_n} \mc  D^{x+z,y+z}(f) \leq \frac{\Lambda_l}{\Lambda_n} \mc D(f) \leq \Big(\frac{2l+1}{2n+1}\Big)^d \frac{K n^{d-\alpha}}{2t }.
\] 

By the convexity of $\mc D^{x,y}(\cdot)$, $\mc D_l(\cdot)$ is also convex and $\mc D_l(\bar f_l) \leq \mc D_l(\bar f)$. We conclude that $\mc D_l(\bar f_l) \leq c(l)/n^\alpha$, where $c(l)$ is a constant that only depends on $l$, $K$ and $t$. Therefore, the integral in (\ref{ec15}) can be estimated by
\begin{equation}
\label{ec16}
\sup_{\mc D_l(f) \leq c(l)/n^\alpha} \int V_0^{l,a}(\xi)  f(\xi) d\nu_\rho.
\end{equation}

Due to the indicator function in the definition of $V_0^{l,a}$, we can restrict the supremum to densities supported on $\{\xi^l(0) <a\}$. This set of densities is compact, and both the function $\mc D_l(f)$ and the  integral in (\ref{ec16}) are continuous. Therefore, the limiit as $n \to \infty$ of (\ref{ec16}) is bounded above by
\begin{equation}
\label{ec17}
\sup_{\mc D_l(f) = 0} \int V_0^{l,a}(\xi) f(\xi) d\nu_\rho.
\end{equation}

If $\mc D_l(f)=0$, then the density $f$ is constant on each ergodic component of the state space $\Omega_l$ with respect to the dynamics generated by $L_{zr}^l$. Since the number of particles is the only conserved quantity with respect to this dynamics, these ergodic components are exactly the sets $\{\xi \in \Omega_l; \xi^l(0) = k\}$, $k \geq 0$. Denote by $\nu_{k,l}(\cdot)$ the measure $\nu_\rho(\cdot|\xi^l(0)=k/\Lambda_l)$. Then, the limit as $l \to \infty$ of the supremum in (\ref{ec17}) is equal to
\[
\lim_{l \to \infty} \sup_{k \leq a \Lambda_l} \int V_0^l(\xi) d\nu_{k,l}.
\]

This last line is exactly the same limit appearing in pg. 89 of \cite{KL}, and it is equal to 0 due to the equivalence of ensembles. This ends the proof of Lemma \ref{1b}.

\subsection{The 2-blocks estimate}
The proof of Lemma \ref{2b} is similar to the proof of the one-block estimate explained in the previous section. When $g(\cdot)$ satisfies {\bf (FEM)} (see section 5.5 of \cite{KL}), using the entropy inequality we can prove that
\begin{multline}
\label{ec18}
\lim_{A \to \infty} \limsup_{\epsilon \to 0} \limsup_{n \to \infty} \sup_{y,z} \int \frac{1}{n^d} \sum_{|x| \leq Mn} \big| \xi^l(x+z) -\xi^l(x+y)\big| \times \\
	\times {\mathbf 1}\{\max\{\xi^l(x+y),\xi^l(x+z)\} \geq a\}  \bar f_t^n(\xi) d\nu_\rho =0.
\end{multline}

When $g(\cdot)$ satisfies {\bf (B)}, we need to use hypothesis {\bf (C)} to introduce this cut-off. The integrand in (\ref{ec18}) is increasing and $\bar f_t^n \nu_\rho \preceq \nu_{\rho'}$.  Therefore, we can estimate the integral replacing $\bar f_t^n(\xi) d\nu_\rho$ by $d\nu_{\rho'}$. By Tchebyshev's inequality, the integral can be estimated by $4 \int \xi(x)^2 d\nu_{\rho'}/a$, which proves that (\ref{ec18}) also holds when $g(\cdot)$ satisfies {\bf (B)}. Define $B_{y,z}^l(\xi) = \mathbf{1}\{\xi^l(z) <a\} \mathbf{1}\{\xi^l(y) <a\}$. As we did in the proof of the one-block estimate, passing the sum to the density $f$, the integral in (\ref{ec12.1}) is equal to
\begin{equation}
\label{ec19}
\Big(\frac{2n+1}{n}\Big)^d \int \big| \xi^l(z) - \xi^l(y)\big| B_{y,z}^l(\xi) \bar f(\xi) \nu_\rho(\xi).
\end{equation}

Notice that the function $|\xi^l(z)-\xi^l(y)|$ depends on the configuration of particles inside the two disjoint blocks $\tau_y \Lambda_l$ and $\tau_z \Lambda_l$. The name ``two-blocks estimate'' comes from this observation. Let us denote by $\Lambda_l^*$ the union of these two blocks, and let us write $\xi_1(x)=\xi(x+y)$, $\xi_2(x) = \xi(x+z)$. We represent $\Lambda_l^*$ as $\Lambda_l \times \Lambda_l$, dropping $y$ and $z$ from the notation. The integral  in (\ref{ec19}) is equal to
\begin{equation}
\label{ec20}
\int_{\Lambda_l \times \Lambda_l} \big| \xi^l_1(0) - \xi_2^l(0)\big|B_l^*(\xi)  \bar f_l^* (\xi_1,\xi_2) \nu_\rho(d\xi_1 \times d\xi_2),
\end{equation}
where $\bar f_l^*= E[\bar f|\mc F_{\Lambda_l^*}]$, $\mc F_{\Lambda_l^*}$ is the $\sigma$-algebra generated by $\{\xi(x); x\in \Lambda_l^*\}$ and $B_l^*(\xi) =   \mathbf{1}\{\xi^l_1(0) <a\} \mathbf{1}\{\xi^l_2(0) <a\}$.

The integral in (\ref{ec20}) depends on $y$, $z$ only through $\bar f_l^*$. The estimation of this integral is now a finite-dimensional problem. This integral is similar to the one in (\ref{ec15}), so we will estimate a suitable version of the Dirichlet form of $\bar f_l^*$. Let us define
\begin{align*}
\mc D_l^*(f) 
	& = \sum_{|x|,|x'| \leq l} \Big(\mc D^{x+y,x'+y}(f) + \mc D^{x+z,x'+z}(f)\Big) + \frac{1}{p(z-y)} \mc D^{y,z}(f) \\
	& =\mc D_{l,0}^*(f) + \frac{1}{p(z-y)} \mc D^{y,z}(f).
\end{align*}

Following the computations made when we estimated $\mc D_l(\bar f_l)$ in the previous section, we see that
\begin{equation}
\label{ec21}
\mc D_{l,0}^*(\bar f_l^*) \leq \frac{c(l)}{n^\alpha}.
\end{equation}

The term $\mc D^{y,z}(\bar f_l^*)$ is the one that connects the behavior of $\bar f_l ^*$ between the two blocks, so it is the most relevant in the definition above. Notice that the dynamics associated to $\mc D_l^*(\cdot)$ corresponds to a system on which particles perform a zero-range process with long jumps restricted to each one of the two boxes, and on which particles can jump from site $y$ to site $z$ (site $z$ to site $y$ resp.) with rate $g(\xi(y))$ ($g(\xi(z))$ resp.) The rate at which particles jump between $y$ and $z$ does not depend on the distance between $y$ and $z$. Now we are ready to state what we mean by the {\em moving particle lemma}:

\begin{lemma}[Moving particle lemma]
\label{l3}
Under the hypothesis of Theorem \ref{t2},
\begin{equation}
\label{ec22}
\lim_{\epsilon \to 0} \limsup_{n \to \infty} \sup_{y,z} \frac{1}{p(z-y)} \mc D^{y,z}(\bar f_l^*) =0,
\end{equation}
where the supremum is over $y, z \in \Lambda_{\epsilon n}$ such that $|z-y| > 2l+1$.
\end{lemma}

If we understand by $p(z-y)^{-1} \mc D^{y,z}(\bar f_l^*)$ as the {\em cost} of moving a particle from $y$ to $z$, what this lemma is saying is that the cost of moving a particle at a macroscopically small distance vanishes when the distance goes to $0$. We will prove this lemma in the next section. Assuming this lemma, it is not difficult to finish the proof of Lemma \ref{2b}. By (\ref{ec21}) and (\ref{ec22}), $\mc D_l^*(\bar f_l^*)$ goes to 0 as $n \to \infty$ and then $\epsilon \to 0$. Therefore, by the same compactness argument used in the proof of the one-block estimate, we can bound the limit of the integral in (\ref{ec20}) by
\[
\lim_{l \to \infty} \sup_{k \leq 2a/\Lambda_l} \int_{\Lambda_l \times \Lambda_l} \big| \xi^l_1(0) - \xi_2^l(0)\big|B_l^*(\xi)   \nu_{k, l}^*(d\xi_1 \times d\xi_2),
\]
where $\nu_{k,l}^*$ is the measure $\nu_\rho$ in $\Lambda_l \times \Lambda_l$, conditioned to have exactly $k$ particles. This limit is equal to 0 due to the equivalence of ensembles, which ends the proof of Lemma \ref{2b}.

\section{The moving particle lemma}
\label{s5}
In this section we prove Lemma \ref{l3}. To avoid heavy notation, we start assuming $d=1$. Later we explain how to generalize the proof to arbitrary dimensions. Assume, without loss of generality, that $y<z$. Let $\{y=y^0, y^1,\dots,y^m=z\}$ be a path from $y$ to $z$. Notice that
\[
\xi^{y,z} = \left( \left( \left( \xi^{y^0,y^1} \right)^{y^1,y^2} \right)\dots \right)^{y^{m-1,m}}.
\]

In other words, moving a particle from $y$ to $z$ is the same that moving a particle from $y$ to $y^1$, then from $y^1$ to $y^2$, etc. Using the formula above, write 
\[
\sqrt{f(\xi^{y,z})} - \sqrt{f(\xi)} = \sum_{i=1}^m \Big\{\sqrt{f((\xi^{y,y^{i-1}})^{y^{i-1},y^i})} - \sqrt{f(\xi^{y,y^{i-1}})}\Big\}.
\]

Using the inequality $(a_1+\dots+a_m)^2 \leq m (a_1^2+\dots+a_m^2)$ and performing a change of variables, we conclude that
\begin{equation}
\label{ec23}
\frac{\mc D^{y,z}(f)}{p(z-y)} \leq m \sum_{i=1}^m \frac{\mc D^{y^{i-1},y^i}(f)}{p(y^{i}-y^{i-1})}
\end{equation}
for any density $f$.
Now the idea is to choose the path between $y$ and $z$ in an adequate way. The simplest choice is $y^i = y+i$. In that case, 
\[
\frac{\mc D^{y,z}(f)}{p(z-y)} \leq |z-y| \sum_{i=1}^{|z-y|} \frac{\mc D^{y+i-1,y+i}(f)}{p(1)}.
\]

Applying this estimate to $\bar f_l^*$, we see that
\begin{align*}
\frac{\mc D^{y,z}(\bar f_l^*)}{p(z-y)}
	&\leq \frac{1}{n} \sum_{|x|\leq n} \frac{\mc D^{x+y,x+z}(f)}{p(z-y)} \leq \frac{|z-y|}{n} \sum_{\substack{|x|\leq n\\i=1,\dots,|z-y|}} \frac{\mc D^{x+y^{i-1},x+y^i}(f)}{p(1)}\\
	& \leq \frac{|z-y|^2}{p(1) n} \mc D(f) \leq \frac{4K\epsilon^2 n^{2-\alpha}}{p(1)}.
\end{align*}

The last but one inequality follows from the fact that each of the terms of the form $\mc D^{x+y^{i-1},x+y^i}(f)$ appears at most $|z-y|$ times in the sum. Notice that in the diffusive case, $\alpha =2$ and this estimate is enough to prove the lemma. The point is that in the previous estimate we have only used jumps of lenght one to move the particle from $y$ to $z$. Therefore, we need to use jumps of various lenghts to get the right estimate. Fix a positive  integer $k$. Assume for a moment that $z-y$ is divisible by $k$. Defining $y^i=y +ik$, (\ref{ec23}) gives us the estimate
\[
\frac{\mc D^{y,z}(\bar f_l^*)}{p(z-y)} \leq 4\epsilon^2 n k^{\alpha-1} \sum_{|x-x'| = k} \mc D^{x,x'}(f).
\]

We make use of the following observation: if $a_k$, $b_k$, $\beta$, $\gamma$ are non-negative numbers such that $a_1+\dots+a_m \leq \gamma$ and $\beta \leq b_k a_k$ for any $k$, then
\[
\beta \leq \frac{\gamma}{\sum_{k=1}^m b_k^{-1}}.
\] 

In our case we take $a_k =\sum_{|x-x'|=k} \mc D^{x,x'}(f)$, $b_k=4\epsilon^2 n k^{\alpha-1}$, $\beta = \mc D^{y,z}(\bar f_l^*)$ and $\gamma = \mc D(f)$. Therefore,
\[
\frac{\mc D^{y,z}(\bar f_l^*)}{p(z-y)} \leq \frac{4K \epsilon^2 n^{2-\alpha}}{\sum_{k=1}^{\epsilon n} k^{1-\alpha}} \leq C(\alpha) \frac{K \epsilon^2 n^{2-\alpha}}{(\epsilon n)^{2-\alpha}} = C(\alpha) K \epsilon^\alpha,
\]
where $C(\alpha)$ is a constant that only depends on $\alpha$. Here we used the fact that $\sum_{k=1}^{\epsilon n} k^{1-\alpha}$ is of the same order than $\int_1^{\epsilon n} x^{1-\alpha}dx$.
Notice that this estimate is still not justified, since we assumed that $k$ is a divisor of $z-y$ for any $k$ between $1$ and $\epsilon  n$. When $k$ is not a divisor of $z-y$, the path from $y$ to $z$ must contain a jump of size different from $k$. The idea is that we can restrict the range of summation of $k$ to an interval of the form $[a\epsilon n, b\epsilon n]$, since in that case the sum will still be of order $(\epsilon n)^{2-\alpha}$. Define $m = |z-y|$. We can assume, without loss of generality, that $m$ is divisible by $6$. In fact, for $l >6$, what we can do is to take $z'$ at distance at most $5$ from $z$ such that $z'-y$ is divisible by 6, and to consider $\mc D^{y,z'}(f)$ instead of $\mc D^{y,z}(f)$.
We write $m= 6m_0$. For $k = 2 m_0 +j$, $j = 1,\dots,m_0$, we have $m = 2k + 2(m_0 -j)$. Therefore, we can move a particle from $y$ to $z$ by making two jumps of length $k$ and two jumps of lenght $m_0 -j$. Notice that $k$ runs from $2 m_0+1$ to $3 m_0$.

In this way we have gained control over the rest of the division of $m$ by $k$. The good property of this decomposition is that for $k \neq k'$ we have $j \neq j'$ and therefore we are not repeating jump lenghts for different $k$'s. Therefore, for any density $f$ we have
\begin{align*}
\frac{\mc D^{y,z}(f)}{p(z-y)} 
	&\leq 4\bigg\{ \sum_{i=1,2} \frac{\mc D^{y^{i-1},y^i}(f)}{p(2m_0 +j)} + \sum_{i=3,4} \frac{\mc D^{y^{i-1},y^i}(f)}{p(m_0 -j)}\bigg\}\\
	&\leq \frac{108 (m_0)^{1+\alpha}}{p(1)} \sum_{i=1,\dots,4} \mc D^{y^{i-1},y^i}(f),
\end{align*}

where $y^0=y$, $y^1=y+k$, $y^2=y+2k$, $y^3= z-(m_0 -j)$ and $y^4=z$. Repeating the calculations we did to estimate $\mc D_l(\bar f_l)$ in the proof of the 1-block estimate, we obtain that
\[
\frac{\mc D^{y,z}(\bar f_l^*)}{p(z-y)} \leq \frac{216 m_0^{1+\alpha}}{p(1)n} \sum_{|x| \leq n(1+\epsilon)} \Big\{ \mc D^{x,x+k}(f) + \mc D^{x,x+m_0 -j}(f)\Big\}.
\]

Since the jumps have different lenghts for different choices of $k$, we can perform a sum from $j=1$ to $j=m_0$ to obtain the bound
\[
\frac{\mc D^{y,z}(\bar f_l^*)}{p(z-y)} \leq  \frac{216K m_0^{\alpha}}{p(1)n^\alpha} \leq \frac{216 K \epsilon^\alpha}{p(1)},
\]
where we have used the estimate $m_0 \leq \epsilon n$ in the last inequality. This ends the proof of Lemma \ref{l3} in dimension $d=1$.

In dimension $d>1$, we proceed as follows. First observe that due to the positivity and continuity of $q(\cdot)$ on the sphere $\bb S^{d-1}$,  the Dirichlet form $\mc D(f)$ is equivalent to the one associated to the zero-range process associated to the transition rate $p(z)=c/\|z\|^{d+\alpha}$. In particular, there is a finite and positive constant $\epsilon_0$ such that $\epsilon_0/\|y-x\| \leq p(y-x) \leq \epsilon_0^{-1}/\|y-x\|^{d+\alpha}$
To simplify the notation, take $d=2$. Assume first that $(z-y)\cdot e_2=0$, where $e_2=(0,1)$ is the second element of the canonical basis in $\bb R^d$. Define $m = |z-y|$. Using the same argument above, we can assume that $m =6m_0$ for some positive integer $m_0$. Define $\mathbf{m_0} = (m_0,0)$ and for $i,j=1,\dots,m_0$, define $\mathbf{k} = 2\mathbf{m_0} + (i,j)$, $\mathbf{k'} = \mathbf{m_0}-(i,j)$. We can move a particle from $y$ to $z$ by making two jumps in direction $\mathbf{k}$ and then two jumps in direction $\mathbf{k'}$. Now we observe that when $i,j$ run from $1$ to $m_0$, the vectors $\mathbf{k}$, $\mathbf{k'}$ are all different. Repeating the computations done in the $d=1$ case,
\[
\frac{\mc D^{y,z}(\bar f_l^*)}{p(z-y)} \leq \frac{C(d,\alpha) m_0^{d+\alpha}}{n^d} \sum_{|x| \leq n(1+\epsilon)} \Big\{ \mc D^{x,x+\mathbf{k}}(f) + \mc D^{x,x+\mathbf{k'}}(f)\Big\}.
\]

Here $C(d,\alpha)$ is a constant which only depends on $\alpha$, the dimension $d$ and $\epsilon_0$. Although we are assuming $d=2$, we include the dependence in dimension of this estimate for clarity. 
Performing a sum over $i,j \in \{1,\dots,m_0\}$, we conclude that
\[
\frac{\mc D^{y,z}(\bar f_l^*)}{p(z-y)} \leq  \frac{C(d,\alpha) K m_0^{\alpha}}{n^\alpha} \leq C(d,\alpha) K \epsilon^\alpha.
\]

When $(z-y)\cdot e_1=0$, we simply exchange the role of the first and second coordinates and we obtain the same estimate. For $y=(y_1,y_2)$, $z=(z_1,z_2)$ in general position,  we define $y' = (z_1,y_2)$. Then, $\mc D^{y,z}(f)/p(z-y) \leq 2\big\{\mc D^{y,y'}(f)/p(y'-y) + \mc D^{y',z}(f)/p(z-y')$, and these two last terms can be estimated as above.

\section{The energy estimate}
\label{s6}
In this section we prove Theorem \ref{t5.1}. 
We start introducing some notation. For a function $F : \bb Z^d \times \bb Z^d \to \bb R$ and a configuration $\xi \in \Omega_{zr}$, we define
\[
\mc E_n^g(\xi,F) =\frac{1}{n^d} \sum_{x,y \in \bb Z^d} n^\alpha p(y-x) (g(\xi(y))-g(\xi(x)))
F_{x,y},
\]
\[
\Gamma_n(\xi,F) = \frac{1}{n^d} \sum_{x,y \in \bb Z^d} n^\alpha p(y-x) \big\{g(\xi(x))+g(\xi(y))\big\} F_{x,y}^2.
\]

The proof of Theorem \ref{t5.1} consists basically on combining the Replacement Lemma \ref{t5} and the following lemma:

\begin{lemma}
\label{l2}
Under the hypothesis of Theorem \ref{t2}, for any finite sequence of functions $F^i_\cdot:  [0,T] \times \bb Z^d \times \bb Z^d \to \bb R$, $i=1,\dots,k$ such that $F^i_{x,y}(t)=0$ if $x,y$ are big enough, we have
\begin{equation}
\label{ec24}
\bb E_{\nu^n} \Big[ \sup_{1\leq i \leq k} \int_0^T \Big\{ \mc E_n^g(\xi_t^n,F^i(t)) - \Gamma_n(\xi_t^n,F^i(t))\Big\} dt \Big] \leq K + \frac{\log k}{n^d}.
\end{equation}
\end{lemma}

\begin{proof}

By the entropy inequality with $\gamma = n^d$, the expectation in (\ref{ec9}) is bounded by
\begin{align*}
K + \frac{1}{ n^d} \log \bb E^{\rho}\Big[\exp\Big\{  n^d \sup_{1\leq i \leq k} \int_0^T \Big\{ \mc E_n^g(\xi_t^n,F^i(t)) - \Gamma_n(\xi_t^n,F^i(t))\Big\} dt \Big\}\Big].
\end{align*}

By estimate (\ref{ec9}), the expression above is bounded by
\[
K+\frac{\log k}{ n^d} +\sup_{1\leq i \leq k} \frac{1}{ n^d} \log \bb E^{\rho}\Big[\exp\Big\{ n^d  \int_0^T \Big\{ \mc E_n^g(\xi_t^n,F^i(t)) - \Gamma_n(\xi_t^n,F^i(t))\Big\} dt \Big\}\Big].
\]

Therefore, it is enough to prove that
\begin{equation}
\label{ec10}
\frac{1}{n^d} \log \bb E^{\rho}\Big[\exp\Big\{ n^d  \int_0^T \Big\{ \mc E_n^g(\xi_t^n,F(t)) - \Gamma_n(\xi_t^n,F(t)\Big\} dt \Big\}\Big] \leq 0
\end{equation}
for any trajectory $F(t)$. Since now the expectation is with respect to the process in equilibrium, powerful variational methods are available to estimate this expectation. By Feynman-Kac formula, the left-hand side of (\ref{ec10}) is bounded by $n^{-d}\int_0^T \lambda(t)dt$, where $\lambda(t)$ is the largest eigenvalue in $\mc L^2(\nu_\rho)$ of the operator $n^\alpha L^{zr} + n^d \mc V_t^n$ and $\mc V_t^n$ is the multiplication operator given by
\[
\mc V_t^n f(\xi) = \Big\{ \mc E_n^g(\xi,F(t)) - \Gamma_n(\xi,F(t))\Big\}f(\xi). 
\]

By the variational formula of the largest eigenvalue of an operator in $\mc L^2(\nu_\rho)$, we see that the left-hand side of (\ref{ec10}) is bounded by
\begin{equation}
\label{ec11}
\int_0^T \sup_{f \in \mc L^2(\nu_\rho)} \big\{ \<\mc V_t^n, f^2\>_\rho -\frac{1}{ n^{d-\alpha}} \<f,-L^{zr} f\>_\rho\big\} dt,
\end{equation}
where $\<\cdot,\cdot\>_\rho$ denotes the inner product in $\mc L^2(\nu_\rho)$. This integral only involves the equilibrium properties of the dynamics. Notice that  $\<f,-L^{zr} f\>_\rho =  \mc D(f^2)$.
Let us estimate each one of the terms in $\<\mc E_n^g(\xi,F(t)),f(\xi)^2\>_\rho$ separatedly.  Define $(\xi+\delta_x)(z) = \xi(z)+ \delta_x(z)$. Using the change of variables $\xi \to \xi + \delta_x$, we see that
\[
\int \big(g(\xi(y))-g(\xi(x))\big) f(\xi)^2 d\nu_\rho = \phi(\rho) \int\big\{f(\xi+\delta_y)^2-f(\xi+\delta_x)^2\big\}d\nu_\rho.
\]

The same change of variables shows that 
\[
\mc D^{x,y}(f^2) = \phi(\rho)\int p(y-x) \big(f(\xi+\delta_y) - f(\xi+\delta_x)\big)^2 d\nu_\rho.
\]

Notice the similarity between these two formulas. Using the formula $a^2-b^2=(a+b)(a-b)$ and the weighted Cauchy-Schwartz inequality $2|ab| \leq \beta a^2 +b^2/\beta$, we see that
\begin{align*}
\big|f(\xi+\delta_y)^2-f(\xi+\delta_x)^2 \big| 
	&= \{f(\xi+\delta_y)+f(\xi+\delta_x)\}\big|f(\xi+\delta_y)-f(\xi+\delta_x)\big| \\
	&\leq \frac{\beta_{xy}}{2} \big\{f(\xi+\delta_y)+f(\xi+\delta_x)\big\}^2\\
	&\quad +\frac{1}{2\beta_{xy}} \big\{f(\xi+\delta_y)-f(\xi+\delta_x)\big\}^2.
\end{align*}

The idea is to choose $\beta_{xy}$ in such a way that the second term above cancels with $n^{\alpha-d} \mc D^{x,y}(f)$. This happens choosing $\beta_{xy} = |F_{xy}(t)|/2$. With this choice for $\beta_{xy}$, we have
\begin{multline*}
\<\mc E_n^g(\xi,F(t)), f^2\>_\rho -\frac{1}{ n^{d-\alpha}} \<f,-L^{zr} f\>_\rho \leq \\
	\leq \frac{1}{n^d}\sum_{x,y  \in \bb Z^d}  n^{\alpha} p(y-x) F_{xy}^2 \int \big\{g(\xi(x))+g(\xi(y))\big\} f(\xi)^2 d\nu_\rho. 
\end{multline*}

This last term is exactly equal to $\<\Gamma_n(\xi,F(t)),f^2\>$, which proves that the supremum in (\ref{ec11}) is less or equal than $0$.
\end{proof}

Notice that the term $\mc E_n^g(\xi,F)$ can be written as
\[
\mc E_n^g(\xi,F) = \frac{1}{n^d} \sum_{x \in \bb Z^d} g(\xi(x)) \sum_{y \in \bb Z^d} n^\alpha p(y-x)(F_{y,x}-F_{x,y}).
\]

For a given function $F$, let us define the symmetric part $F^s$ and antisymmetric part $F^a$ of $F$ by
\[
F^s_{x,y} =\frac{1}{2} \big(F_{x,y} + F_{y,x}\big),
\]
\[
F^a_{x,y}=\frac{1}{2} \big(F_{x,y}-F_{y,x}\big).
\]

We have  $\mc E_n^g(\xi,F) = \mc E_n^g(\xi,F^a)$ and $\Gamma_n(\xi, F) \geq \Gamma_n(\xi,F^a)$. Therefore, the estimate (\ref{ec24}) is better for antisymmetric functions $F^i$. For an antisymmetric function $F$, we have
\[
\mc E_n^g(\xi,F) = \frac{1}{n^d} \sum_{x \in \bb Z^d} g(\xi(x)) \sum_{z \in \bb Z^d} n^\alpha p(z)(F_{x+z,x}+F_{x-z,x}).
\]

Let us denote by $\mc C_{b,\text{ant}}^2(\bb R^d)$ the set of uniformly continuous, bounded, twice differentiable functions $G: \bb R^d \times \bb R^d \to \bb R^d$ such that $G(x,y) = -G(y,x)$. Notice that in particular $G(x,x)=0$ for any $G \in \mc C_{b,\text{ant}}^2(\bb R^d)$.
Let $\{G^i;i=1,\dots,k\}$ be a sequence of continuous paths $G^i: [0,T] \to \mc C_{b,\text{ant}}^2(\bb R^d)$ and define $F^i_t(x,y) = G_t^i(x/n,y/n)$. 
We will extend the definitions of $\mc L_n G$, $\mc L G$ to antisymmetric functions. For $G \in \mc C^2_{b,\text{ant}}(\bb R^d)$ we define
\[
\mc L_n G^i_t(x/n) = \frac{1}{n^d} \sum_{z \in \bb Z^d} q(z/n) \Big\{G^i_t\Big(\frac{x+z}{n},\frac{x}{n}\Big) + G^i_t\Big(\frac{x-z}{n},\frac{x}{n}\Big) \Big\},
\]
\[
\mc L G(x) = \int q(z) \big\{G(x+z,x)+G(x-z,x)\big\}dz.
\]

With this notation, 
\[
\mc E_n^g(\xi,F_t^i) = \frac{1}{n^d} \sum_{x \in \bb Z^d} g(\xi(x)) \mc L_n G^i_t(x/n).
\]

Using the second-order Taylor expansion of $G^i_t$ around $(x/n,x/n)$, we see that 
\[
\lim_{n \to \infty} \sup_{x \in \bb Z^d} \big| \mc L_n G_t^i(x/n) -\mc L G_t^i(x/n)\big| =0.
\]

These definitions reduce to the previous definitions of $\mc L_n$, $\mc L$ when we take $G^i_t(x,y) = \mc G^i_t(y)-\mc G^i_t(x)$.
In the same way, we define
\[
\mc Q_n G^i_t(x/n) = \frac{1}{n^d} \sum_{z \in \bb Z^d} q(z/n) G((x+z)/n,x/n)^2,
\]
\[
\mc Q G(x) = \int q(z) G(x+z)^2 dz.
\]

Then,
\[
\Gamma_n(\xi,F^i_t) = \frac{2}{n^d} \sum_{x \in \bb Z^d} g(\xi(x)) \mc Q_n G^i_t(x/n)
\]
and  $\mc Q_n G^i_t(x/n)$ is converging uniformly to  $\mc Q G^i_t(x/n)$ in the sense that 
\[
\lim_{n \to \infty} \sup_{x \in \bb Z^d} \big| \mc Q_n G_t^i(x/n) -\mc Q G_t^i(x/n)\big| =0.
\]

Let the sequence $\{G^i_t,i=1,\dots,k\}$ be fixed. Assume that the trajectories $G^i_t$ are of class $\mc C^1$. As in Section \ref{s3.2}, let $\pi_\cdot$ be a limit point of the empirical process, and let $n'$ be a subsequence such that $\pi_\cdot^{n'}$ converges to $\pi_\cdot$. Remember that we already proved that $\pi_\cdot(dx) = u(x,\cdot) dx$ for some locally integrable density $u(x,t)$. Using the Replacement Lemma in the same way we did it in Section \ref{s3.2}, we can prove that
\begin{multline*}
\lim_{n' \to \infty} \sup_{1\leq i\leq k} \int_0^T \Big\{ \mc E_{n'}^g(\xi_t^{n'},F^i(t)) - \Gamma_{n'}(\xi_t^{n'},F^i(t)) \Big\} dt \\
	= \sup_{1\leq i\leq k} \int_0^T \int \phi(u(x,t)) \big\{\mc L G^i_t(x) - 2\mc Q G^i_t(x)\big\}dx dt.
\end{multline*}

Combining this last line with the energy estimate (\ref{ec24}), we conclude that 
\[
E\Big[ \sup_{1\leq i\leq k} \int_0^T \int \phi(u(x,t)) \big\{\mc L G^i_t(x) - 2\mc Q G^i_t(x)\big\}dx dt \Big]
\leq K.
\]

By the monotone convergence theorem, we conclude that
\[
E\Big[ \sup_{G} \int_0^T \int \phi(u(x,t)) \big\{\mc L G_t(x) - 2\mc Q G_t(x)\big\}dx dt \Big]
\leq K,
\]
where the supremum is over the set of continuous paths $G:[0,T] \to \mc C_{b,\text{ant}}^2(\bb R^d)$. We conclude that the limit density $u(x,t)$ satisfies 
\begin{equation}
\label{ec25}
\sup_{G} \int_0^T \int \phi(u(x,t)) \big\{\mc L G_t(x) - 2\mc Q G_t(x)\big\}dx dt <+\infty,
\end{equation}
almost surely with respect to the law of $\pi_\cdot$. Let us define $\phi_t(x) = \phi(u(x,t))$. Let $\|\cdot\|_{\phi,T}$ be the norm defined by
\begin{align*}
\| G \|_{\phi,T}^2 
	&= 2\int_0^T \int \phi_t(x) \mc Q G_t(x) dx dt\\
	&= \int_0^T \iint \big(\phi_t(x)+\phi_t(y)\big) h(y-x) G_t(y,x)^2 dx dy dt
\end{align*}
for any continuous trajectory $G: [0,T] \to \mc C_{b,\text{ant}}(\bb R^d)$. The second definition will be more convenient by symmetry reasons. We denote by $\mc L^2_{\phi,T,\text{ant}}$ the Hilbert space obained as the closure of such trajectories under $\| \cdot \|_{\phi,T}$. With this notation, we can rewrite (\ref{ec25}) as
\[
\sup_{G} \Big\{\int_0^T \int \phi(u(x,t)) \mc L G_t(x) dx dt - \|G\|_{\phi,T}^2\Big\} <+\infty.
\]

We recognize in this formula the variational formula for the norm with respect to $\mc L^2_{\phi,T,\text{ant}}$ of the linear funcional $\lambda(G)=1/2 \int_0^T \int \phi(u(x,t)) \mc L G_t(x) dx dt$. We conclude that $\lambda(\cdot)$ is a continuous functional in $\mc L^2_{\phi,T,\text{ant}}$. By Riesz representation teorem, there exists an antisymmetric function $F_t(x,y)$ such that
\[
\int_0^T \iint \big(\phi_t(x)+\phi_t(y)\big) h(y-x) F_t(x,y)^2 dx dt <+\infty, \text{ and}
\]
\[
\lambda(G) = \frac{1}{4}\Big\{ \|F+G\|_{\phi,T}^2-\|F-G\|_{\phi,T}^2\Big\}.
\]

The arguments above allow to justify the formal computations of calculus of variations, from which we obtain that
\[
F_t(x,y) = \frac{1}{2} \frac{\phi_t(y)-\phi_t(x)}{\phi_t(y)+\phi_t(x)}, \text{ and therefore}
\]
\[
	\int_0^T \iint h(y-x) \frac{\big(\phi_t(y)-\phi_t(x)\big)^2}{\phi_t(x)+\phi_t(y)} dx dy dt <+\infty.
\]

Aside from the factor $\phi_t(y)+\phi_t(x)$ in the denominator, this estimate is exactly the estimate in Theorem \ref{t5.1}. At this point we need an extra argument. When $g(\cdot)$ satisfies {\bf (B)}, the interaction rate $\phi$ is bounded above by $\phi_c$ and therefore
\[
\int_0^T \mc E(\phi(u_t),\phi(u_t)) dt \leq 2 \phi_c \int_0^T \iint h(y-x) \frac{\big(\phi_t(y)-\phi_t(x)\big)^2}{\phi_t(x)+\phi_t(y)} dx dy dt <+\infty.
\]

When $g(\cdot)$ satisfies {\bf (FEM)}, $\phi_c=+\infty$ and another argument is needed. We appeal to Hypothesis {\bf (C)}. We point out that when $g(\cdot)$ satisfies {\bf (FEM)}, this is the only place where we need to consider an increasing rate function $g(\cdot)$. In Appendix \ref{B} we explain how to get rid of Hypothesis {\bf (C)}, but still assuming that $g(\cdot)$ is increasing. Under Hypothesis {\bf (C)}, for any function $G \in \mc C_c(\bb R^d)$ we have
\[
\limsup_{n \to \infty} \pi_t^n(G) \leq \limsup_{n \to \infty} \frac{1}{n^d} \sum_{z \in \bb Z^d} \rho' G(z/n) = \rho' \int G(x) dx,
\]
and in particular $u(x,t) \leq \rho'$ for any $x \in \bb R^d$, $t \in [0,\infty)$. In this case we have
\[
\int_0^T \mc E(\phi(u_t),\phi(u_t)) dt \leq 2 \phi(\rho') \int_0^T \iint h(y-x) \frac{\big(\phi_t(y)-\phi_t(x)\big)^2}{\phi_t(x)+\phi_t(y)} dx dy dt <+\infty,
\]

This ends the proof of Theorem \ref{t5.1} under both conditions {\bf (FEM)} or {\bf (B)}. For the sake of completeness, now we state a weaker form of Theorem \ref{t2}, which follows after checking which hypothesis we have used on each step of the proof of Theorem \ref{t2}.

\begin{theorem}
\label{t6}
Let $g(\cdot)$ be an interaction rate satisfying $\sup_n|g(n+1)-g(n)| <+\infty$ and {\bf (FEM)}. Let $p(\cdot)$ be a transition rate satisfying {\bf (P)}. Let $u_0: \bb R^d \to [0,\infty)$ be a measurable, locally integrable initial profile and let $\{\nu^n\}_n$ be a sequence of probability measures in $\Omega_{zr}$ associated to $u_0$. Let $\xi_t^n$ be the zero-range process $\xi_{t n^\alpha}$ starting from $\nu^n$. Assume that there are positive, finite constants $\rho, K$ such that
\begin{itemize}
\item[{\bf (H)}] For any $n \geq 0$, $H(\nu^n |\nu_\rho) \leq K n^d$. 
\end{itemize}

Let $\pi_\cdot^n$ be the empirical measure associated to $\xi_t^n$. Then, the sequence $\{\pi_\cdot^n\}_n$ is tight, and any limit point $\pi_\cdot$ of $\{\pi_\cdot^n\}_n$ is concentrated on paths of the form $u(x,t)dx$, where $u(x,t)$ is a weak solution of (\ref{ec2}) satisfying \eqref{ec2.4}, \eqref{ec2.2} and in place of \eqref{ec2.3},
\begin{equation}
\label{ec25.1}
\int_0^T \mc E(\sqrt{\phi(u(\cdot,t))},\sqrt{\phi(u(\cdot,t))}) dt <+ \infty
\end{equation}
for any $T >0$.
\end{theorem}

We call these solutions {\em finite entropy} solutions, by analogy with the solutions of heat equation with finite entropy. 
This Theorem implies the hydrodynamic limit as stated in Theorem \ref{t2}, conditioned on a uniqueness result for weak solutions of \eqref{ec2} under \eqref{ec25.1}.

\section{Uniqueness results for the hydrodynamic equation}
\label{s7}

\subsection{The linear case}
\label{s7.1}
In this section we prove the uniqueness results we need in order to establish the hydrodynamic limits of the exclusion process and zero-range process with long jumps. We start with the linear case. We learned this proof from Luis Silvestre. Let $u(x,t)$ be a weak solution of \eqref{ec1}. By linearity, we can assume $u_0(\cdot) \equiv 0$. We will extend $u(x,t)$ to negative values of $t$ by taking $u(x,t)=0$ if $t<0$. For $\tau>0$, we define $\theta_\tau u(x,t)=u(x,t+\tau)$. The function $\theta_\tau u$ is also a weak solution of \eqref{ec1}. Since linear, convex combinations of solutions are also solutions with the same initial condition, for any integrable function $h:[0,a] \to [0,\infty)$,
\[
u^h(x,t)= \int_0^a \theta_\tau u(x,t) h(\tau) d\tau = \int_{-\infty}^\infty u(x,t-s) h(s) ds
\]
is a weak solution of \eqref{ec1}. In particular, if $h$ is the restriction to $[0,a]$ of a continuously differentiable function of compact support contained in $(0,a)$, then $u^h(x,t)$ is continuously differentiable with respect to time. In the same spirit, let $h': \bb R^d \to [0,\infty)$ be a twice continuously differentiable function of compact support. The function
\[
u^{h,h'}(x,t) =\int u^h(x-y,t) h'(y) dy
\]
is also a weak solution of \eqref{ec1} with initial profile $u_0 \equiv 0$. But now the function $u^{h,h'}(x,t)$ is also twice continuously differentiable in time. Let us suppose that $u(x,t)$ is not identically equal to 0. Then, considering $h$, $h'$ as properly defined approximations of the unity, we can assume tha $u^{h,h'}(x,t)$ is also not identically equal to 0. We say that a solution of \eqref{ec1}  is {\em classical} if the solution is twice differentiable in space and differentiable in time. In that case, \eqref{ec1} is satisfied for each pair $(x,t)$, and not only on a weak sense. We have proved that the existence of a weak solution of \eqref{ec1}, not identically null and with initial condition $u_0 \equiv 0$, imply the existence of a classical solution with the same properties.  But for classical, bounded solutions, uniqueness is an immediate consequence of the maximum principle, as we will see. 

Let $u(x,t)$ be the classical solution constructed above. Fix $T>0$ and define $U_T= \sup\{u(x,t);  x \in \bb R^d, t \in [0,T]\}$. If the supremum $U_T$ is attained at some point $(x_0,t_0)$, we necessarily have $\mc L u(x_0,t_0) \leq 0$, with strict inequality if $u(\cdot,t_0)$ is not identically constant. In that case, since $(\partial_t -\mc L)u=0$, we conclude that $\partial_t u(x_0,t_0) <0$ and $u(x_0,\cdot)$ is decreasing in a neighborhood of $t_0$. This contradicts the fact that $(x_0,t_0)$ is the global maximum of $u$, unless that $t_0=0$. This is basically what the maximum principle says. Of course, if $u(\cdot,t_0)$ is constant or if $U_T$ is not attained at any point, we need an extra argument. Let $\epsilon \in (0,1/2T)$ be fixed, and take $(x_0,t_0)$ such that $u(x_0,t_0) \geq U_T-\epsilon$. Consider the test function $g(x) = e^{-x^2/2}$ (actually, any smooth function with a strict maximum at $x=0$ and a fast decay at infinity would be as good as $g(x)$). Consider $\lambda >0$ such that $f(x) =: g(\lambda x)$ satisfies $\sup_x|\mc L f(x)| <\epsilon$ and define $v(x,t) = 2\epsilon f(x-x_0) -2\epsilon^2 t$. Then, $(\partial_t -\mc L)v(x,t) = -2\epsilon^2 +2 \epsilon \mc L f(x-x_0) <0$ for any $x \in \bb R^d$, $t \in [0,T]$. Therefore, $(\partial_t-\mc L)(u+v)= (\partial_t -\mc L)u <0$ for any $x \in \bb R^d$, $tÊ\in [0,T]$. Moreover,
\[
(u+v)(x_0,t_0) > U_T-\epsilon + 2\epsilon -2\epsilon^2 t_0 = U_T+\epsilon -2\epsilon^2t_0\geq U_T, \text{ and }
\]
\[
(u+v)(x,t)<U_T-\epsilon^2 t
\]
for $x$ far enough from $x_0$. Now we can conclude that $u+v$ has a global maximum. Since $(\partial_t -\mc L)(u+v) <0$, the argument exposed above tells us that the global maximum of $u+v$ is attained at $t=0$. But $u(x,0)+v(x,0)= 2\epsilon f(x-x_0) \leq 2\epsilon$, so $u(x,t) \leq 2\epsilon -v(x,t)$ for any $x \in \bb R^d$, $t\in [0,T]$. Since $v(x,t) \leq 2\epsilon -2\epsilon^2t$, we conclude that $u(x,t) \leq 2\epsilon(1+\epsilon T)$ for any $x \in \bb R^d$, $t\in [0,T]$. Since $\epsilon$ is arbitrary, we conclude that $u(x,t) \leq 0$ for any $x,t$. Repeating the argument for $-u(x,t)$, we conclude that $u(x,t) \equiv 0$.

\subsection{The nonlinear case}

In the nonlinear case, uniqueness of solutions can be obtained by an argument due to Oleinik. In order to use Oleinik's argument, we first need to prove that solutions have a bounded second moment with respect to the reference density $\rho$. We start with some elementary lemmas. Here and below, $c$ will denote a constant which may change from line to line, and that depends only on fixed parameters, like $\epsilon_0$, $\alpha$, $d$, etc. For simplicity we will assume the ellipticity condition $\epsilon_1 \leq \phi'(u)$ for any $u \geq 0$, although the arguments can be carried out for any function $\phi$ arising from the interaction rates $g(\cdot)$ considered in this article.

\begin{lemma}
\label{l5}
There is a finite constant $c=c(d,\alpha)$ such that for any $b: \bb R^d \to \bb R$ with $\|b\|_\infty, \|\nabla b\|_\infty<+\infty$, 
\[
\sup_{x \in \bb R^d} \int h(y-x)\big(b(y)-b(x)\big)^2dy \leq c \|b\|_\infty^{2-\alpha}\|\nabla b\|_\infty^{\alpha}.
\]
\end{lemma}

\begin{proof}
This is an easy application of the mean value theorem. In one hand, $(b(y)-b(x))^2 \leq \|\nabla b\|_\infty^2|y-x|^2$. In the other hand, $(b(y)-b(x))^2 \leq 4\| b\|_\infty^2$. The first estimate is good when $|y-x|$ is small. The second one is good when $|y-x|$ is big. After changing to polar coordinates, we have
\begin{align*}
 \int h(y-x)\big(b(y)-b(x)\big)^2dy 
 	&\leq \frac{A_d}{\epsilon_0} \int_0^l \|\nabla b\|_\infty^2r^{1-\alpha} dr + \frac{4 A_d}{ \epsilon_0} \int_l^\infty  \|b\|_\infty^2 r^{-1-\alpha} dr \\
 	&\leq \frac{ A_d \|\nabla b\|_\infty^2 l^{2-\alpha}}{\epsilon_0 (2-\alpha)} 
		+ \frac{4  A_d  \|b\|_\infty^2}{\epsilon_0\alpha l^\alpha}.
\end{align*}
where $A_d$ is the area of the unit sphere in $\bb R^d$ and $l>0$ is arbitrary. Taking $l=\|b\|_\infty/\|\nabla b\|_\infty$, we obtain the desired bound.
\end{proof}

\begin{lemma}[Poincar\'e inequality]
\label{l6}
There exists a positive constant $c$ such that for any locally integrable function $b:\bb R^d \to \bb R$ such that $\mc E(b,b)<+\infty$ and any $M>0$,
\[
\int\limits_{|x|\leq M} b(x)^2 dx \leq cM^\alpha \mc E(b,b) + \frac{1}{M^d} \Big(\int\limits_{|x|\leq M} b(x) dx\Big)^2.
\]
\end{lemma}

\begin{proof}
In our setting, this inequality is specially elementary. In fact, it is enough to observe that
\[
\iint\limits_{\substack{|x| \leq M\\ |y| \leq M}}\big(b(x) -b(y)\big)^2 dx dy \leq \epsilon_0^{-1} M^{d+\alpha} \iint\limits_{\substack{|x| \leq M\\ |y| \leq M}}h(y-x) \big(b(y)-b(x)\big)^2 dx dy,
\]
and to develop the square in the first integral.
\end{proof}

Remember that our objective is to prove that energy solutions of \eqref{ec2} have bounded second moment with respect to the reference density $\rho$. The idea is to prove that suitable truncations have uniformly bounded second moments and then to pass to the limit.  
Let us define $a:\bb R^d \to \bb R$ by
\[
a(x) =
\begin{cases}
1, & |x|\leq1\\
2-|x|,& 1\leq |x| \leq 2\\
0,& 2\leq |x|\\
\end{cases}
\]
and for $M >0$ define $a^M(x)=a(x/M)$. For a given function $u:\bb R^d \to \bb R$, we define its truncation $u^M$ by taking $u^M(x) = a^M(x) u(x)$ for any $x \in \bb R^d$. Let us write $B(M)=\{ x \in \bb R^d: |x|\leq M\}$. Notice that $u^M$ has a support contained in $B(2M)$.

\begin{lemma}
\label{l7}
Let $u:\bb R^d \to \bb R$ be locally integrable and such that $\mc E(u,u) <+\infty$. Then, there is a constant $c$ independent of $u$ such that for any $M >0$, $\mc E(u^M,u^M)$ is finite and
\begin{equation}
\label{ec25.5}
\mc E(u^M,u^M) \leq c\mc E(u,u) +\frac{c}{M^{d+\alpha}}\Big(\int\limits_{B(2M)} u(x) dx \Big)^2.
\end{equation}
\end{lemma} 

\begin{proof}
Let us compute $\mc E(u^M,u^M)$. Since $u^M(x)=0$ if $x \notin B(2M)$ and since the kernel $h(\cdot)$ is symmetric, we see that
\begin{align*}
\mc E(u^M,u^M)
	&= \iint h(y-x)\big(u(x)a^M(x)-u(y)a^M(y)\big)^2dxdy\\
	&\leq 2\iint\limits_{\substack{x \in B(2M)}} h(y-x) \big(u(x)a^M(x)-u(y)a^M(y)\big)^2dxdy. \\
\end{align*}

Let us write $u(x)a^M(x)-u(y)a^M(y)=u(x)\big(a^M(x)-a^M(y)\big)+a^M(y)\big(u(x)-u(y)\big)$. Using the elementary inequality $(a+b)^2 \leq 2a^2+2b^2$, we obtain the bound
\begin{equation}
\label{ec25.6}
\begin{split}
\mc E(u^M,u^M)
	&\leq 4 \iint\limits_{\substack{x \in B(2M) }} h(y-x) u(x)^2\big(a^M(y)-a^M(x)\big)^2 dxdy\\
	&+ 4 \iint\limits_{\substack{x \in B(2M) }} h(y-x) a^M(y)^2\big(u(y)-u(x)\big)^2 dxdy.
\end{split}
\end{equation}

Notice that $\nabla a^M(x) = \nabla a(x/M)/M$. By Lemma \ref{l5}, the first integral is bounded by $cM^{-\alpha} \int_{B(2M)} u(x)^2dx$.
Using the fact that $a^M$ is bounded by $1$, the second integral in \eqref{ec25.6} is bounded by $\mc E(u,u)$. Therefore, we conclude that
\begin{equation}
\label{ec25.7}
\mc E(u^M,u^M) \leq \frac{c}{M^\alpha} \int\limits_{B(2M)} u(x)^2dx + 4\mc E(u,u).
\end{equation}

The lemma follows using Lemma \ref{l6} to estimate the integral in \eqref{ec25.7}.
\end{proof}

Let $u(x,t)$ be an energy solution of \eqref{ec2}. Let us write $u_t(\cdot)=u(\cdot,t)$ and $\phi_t(\cdot) = \phi(\cdot,t)$. Let us denote by $\<\cdot,\cdot\>$ the inner product in $\mc L^2(\bb R^d)$. 
Thanks to the energy estimate \eqref{ec2.3}, relation \eqref{ec2.2} can be written as
\begin{equation}
\label{ec25.2}
\int_0^T \<u_t,\partial G_t\> dt -\int_0^T \mc E(\phi_t,G_t) dt + \<u_0,G_0\> =0.
\end{equation}

When the test function $G_t$ is of the form $G(x,t) = \int_t^T g(x,s)ds$ for a smooth function $g_t$, this relation can be written as
\begin{equation}
\label{ec25.3}
\int_0^T \<u_t,g_t\> dt + \frac{1}{2} \mc E \Big(\int_0^T \phi_t dt, \int_0^T g_t dt\Big) = \int_0^T \<u_0, g_t\> dt.
\end{equation}

Let us assume for a moment that we can choose $g_t = \phi_t$ as a test function. In that case the first term in \eqref{ec25.3} will be comparable to the  $\mc L^2$ norm of $u_t$, by convexity the second term will be bounded by the energy of $u_t$ and the third term will be bounded for suitable initial data $u_0$. Of course, $\phi_t$ is not smooth, so we need to approximate $\phi_t$ in a suitable way. 
The first thing to do is to truncate $\phi_t$.
Let us define $\phi_t^M = a^M \phi_t$ and let $\{i_\epsilon;\epsilon>0\}$ be an approximation of the identity. The function $\phi_t^{M,\epsilon}= i_\epsilon \ast \phi_t^M$ is a legitimate test function, so all the integrals in \eqref{ec25.3} are well defined. Therefore, we have
\begin{equation}
\label{ec27}
\int_0^T \<u_t,\phi_t^{M,\epsilon}\> dt + \frac{1}{2} \mc E \Big(\int_0^T i_\epsilon \ast \phi_t dt, \int_0^T \phi_t^M dt\Big) = \int_0^T \<i_\epsilon \ast u_0, \phi^M_t\> dt.
\end{equation}

Notice that there is a constant $c>0$ such that $u \leq c(\mc H(u)+1)$ for any $u \geq0$. By the entropy bound, we have
\[
\int_0^T \int\limits_{B(2M)} \phi_t(x) dx dt \leq cT(1+M^d).
\]

Therefore, by Lemma \ref{l7} and the convexity of the energy form, we have
\begin{align*}
\mc E\Big(\int_0^T i_\epsilon \ast \phi_t dt, \int_0^T \phi_t^M dt\Big)
	&\leq \frac{1}{2}\Big\{\mc E\Big(\int_0^T i_\epsilon \ast \phi_t dt, \int_0^T i_\epsilon \ast\phi_t dt\Big) \\
	&\quad+ \mc E\Big(\int_0^T \phi_t^M dt, \int_0^T \phi_t^M dt\Big)\Big\}\\
	&\leq c\int_0^T \mc E(\phi_t,\phi_t)dt + \frac{c}{M^{d+\alpha}}\Big(\int_0^T \int\limits_{B(2M)} \phi_t(x) dxdt\Big)^2\\
	&\leq c(1+M^{d-\alpha}).
\end{align*}

For the term involving the initial condition, we have
\[
\int_0^T \<i_\epsilon u_0, \phi_t^M\> dt \leq c \|u_0\|_\infty(1+M^d).
\]

Therefore,  for $M \geq 1$ we have $\int_0^T \<u_t,\phi_t^{M,\epsilon}\> dt  \leq c(1+M^d)$.
By Fatou's lemma, $\int_0^T \<u_t,\phi_t^{M}\> dt \leq\liminf_\epsilon \int_0^T \<u_t,\phi_t^{M,\epsilon}\> dt $. We conclude that $\int_0^T \< u_t, \phi_t^M\>dt \leq c(1+M^d)$. Notice that $\phi_t^M \geq \epsilon_0 u_t$ in $B(M)$. Therefore,
\begin{equation}
\label{ec27.1}
\int_0^T \int\limits_{B(M)} u(x,t)^2 dx dt \leq c(1+M^d)
\end{equation}
 and in particular $u(x,t)$ belongs to $\mc L^2_{\text{loc}}(\bb R^d \times [0,T])$. 
Let us define  $\bar u_t = u_t -\rho$, $\bar \phi_t = \phi_t -\phi(\rho)$. Relation \eqref{ec25.3} can be written as
\begin{equation}
\label{ec27.2}
\int_0^T \<\bar u_t,g_t\> dt + \frac{1}{2} \mc E \Big(\int_0^T \bar \phi_t dt, \int_0^T g_t dt\Big) = \int_0^T \<\bar u_0, g_t\> dt.
\end{equation}

Now we take as a test function $\bar \phi_t^{M,\epsilon} = i_\epsilon \ast (\bar \phi_t a^M)$. We will repeat the computations above. Remember that $\mc H(u) \sim c(u-\rho)^2$ for $u$ close to $\rho$ and $\mc H(u) \sim c u$ for $u$ big. Therefore, there exists a finite constant $c$ such that $\mc H(u) \geq c \min\{|u-\rho|,(u-\rho)^2\}$ for any $u \geq 0$. Therefore, there is a constant $c$ which only depends on $T$, on the entropy bound and on $\int (u_0-\rho)^2dx$, $\|u_0\|_\infty$ such that
$\int_0^T \<\bar u_t,\bar \phi_t^{M,\epsilon}\> dt  \leq c$. This estimate is better than the previous one, due to the introduction of the reference density in the definition of $\bar \phi_t$. 

This time we will use \eqref{ec25.7} instead of Poincar\'e inequality in order to estimate the energy form. We have
\[
\mc E\Big(\int_0^T \phi_t dt, \int_0^T \bar \phi_t^{M,\epsilon}  dt\Big)
	\leq c \int_0^T \mc E(\phi_t,\phi_t) dt + \frac{c}{M^\alpha} \int_0^T \int\limits_{|x|\leq 2M} u(x)^2 dx dt.
\]

The energy bound says that we can bound the first integral by a constant $c$ independent of $\epsilon$ or $M$. Notice that the second bound now makes sense since  $u^2$ is locally integrable. Also due to this fact, we have $\lim_\epsilon \int_0^T \<\bar u_t, \bar \phi_t^{M,\epsilon}\>dt = \int_0^T \<\bar u_t, \bar \phi_t^M\>dt$. Since $\bar u_t$ and $\bar \phi_t^M$ have the same sign, we conclude that
\begin{equation}
\label{ec29}
\int_0^T \int \limits_{B(M)} \bar u_t(x)^2 dx dt \leq c_0 + \frac{c_1}{M^\alpha} \int_0^T \int\limits_{B(2M)} \bar u_t(x)^2 dx dt.
\end{equation}

We will see that inequality \eqref{ec29} plus the entropy estimate \eqref{ec2.4} imply that $\int_0^T \int \bar u_t(x)^2 dx dt <+\infty$. Let us define $a_n = \int_0^T \int_{B(2^n)}\bar u_t(x)^2$. By \eqref{ec29} we have $a_{n+1} \geq (a_n-c_0)2^{\alpha n}/c_1$. The idea is that a sequence $\{a_n\}_n$ satisfying this inequality is either bounded or has a very fast growth. Let us assume that $\lim_n a_n =+\infty$. Let $n_0$ be such that $n_0 \geq (\log_2 c_1 +1)/\alpha$ and $a_{n_0} \geq 2c_0$. If $n \geq n_0$, then $(a_n-c_0)2^{\alpha n}/c_1 \geq (a_n-c_0)2^{\alpha n_0}/c_1 \geq 2(a_n-c_0) = a_n +(a_n-2c_0)$. Therefore, if also $a_n \geq 2c_0$, we conclude that $a_{n+1} \geq a_n$. Inductively we conclude that $a_{n+1}\geq a_n \geq 2c_0$ for any $n \geq n_0$.
Now we start to iterate the inequality in order to get better estimates. 
Since $a_n \geq 2c_0$ for $n \geq n_0$, we have $a_{n+1} \geq 2^{\alpha n} c_0/c_1$. 
Therefore, there exists a constant $\beta_1>0$ such that $a_n \geq \beta_1 2^{\alpha n}+c_0$ for any $n \geq n_0$. Inductively, let us assume that there is a constant $\beta_p>0$ such that $a_n \geq \beta_p 2^{p\alpha n}+c_0$ for any $n \geq n_0$. Then, for $n \geq n_0$ we have $a_{n+1} \geq \beta_p 2^{(p+1)\alpha n}/c_1$ and there exists $\beta_{n+1} >0$ such that $a_n \geq \beta_{p+1} 2^{(p+1)\alpha n} +c_0$ for any $n \geq n_0$. 

We have just proved that if $a_n$ converges to $+\infty$ as $n \to \infty$, then for any $p>0$ there exists a constant $\beta_p$ such that $a_n \geq \beta_p 2^{p\alpha n}$ for any $n$ large enough. In words, the integral of $\bar u_t(x)^2$ over a ball of radius $M$ grows faster than any power of $M$. But this can not happen. In fact, by \eqref{ec27.1}, $a_n \leq c(1+2^{nd})$. Taking $p$ such that $\alpha p>d$ and $n$ big enough, we obtain a contradiction.

We have now proved that any energy solution of \eqref{ec2} is square integrable with respect to the reference density. Let us take two solutions $u^1$, $u^2$ of \eqref{ec2}. Let us write $\phi_t^i(x) = \phi(u^i(x,t))$ for $i=1,2$. By \eqref{ec27.2}, for any test function $g_t$ we have
\[
\int_0^T \< u_t^1-u_t^2,g_t\> dt + \frac{1}{2} \mc E \Big(\int_0^T  (\phi_t^1-\phi_t^2) dt, \int_0^T g_t dt\Big) = 0.
\]

But now we know that $u^1-u^2$ is in $\mc L^2(\bb R^d \times [0,T])$. Therefore, we can approximate $u_t^1-u_t^2$ by test functions to obtain that
\[
\int_0^T \int  (u_t^1(x)-u^2_t(x))^2 dx dt + \frac{1}{2} \mc E \Big(\int_0^T  (\phi_t^1-\phi_t^2) dt, \int_0^T (\phi_t^1-\phi_t^2) dt\Big) = 0.
\]

Since both terms above are non-negative, we conclude that $u^1=u^2$, and uniqueness follows.

\section{The tagged particle problem}
\label{s8}
As an application of the results proved above, we obtain in this section the scaling limit  of a tagged particle in the zero-range process with long jumps. A similar result has been obtained for the exclusion process with long jumps \cite{Jar}  and for a mean zero, finite-range zero-range process in dimension $d=1$ \cite{JLS}.

Take a sequence of initial measures $\{\nu^n\}$ satisfying the hypothesis of Theorem \ref{t2}. Assume as well that $\nu^n(\xi(0) \geq 1)=1$ for any $n$. This condition guarantees that the process $\xi_t^n$ starts with at least one particle at the origin. Tag one of these particles, and follow its evolution, as well as the evolution of the particles as a whole. We need to decide in which way the tagged particle interacts with the other particles. If we do not want the tagged particle to be different from other particles, each time a particle jumps from the site $x$ where the tagged particle is, we decide that the tagged particle is the one who jumps with probability $1/\xi_t^n(x)$.

We obtain in this way a process $(\xi_t^n,X_t^n)$, where $X_t^n$ is the position of the tagged particle at time $t$. We call $X_t^n$ the {\em tagged particle process}. Notice that due to the interaction with other particles, $X_t^n$ is not a Markovian process. However, it can be shown that $X_t^n$ is a local martingale (it is a martingale for $\alpha>1$). We want to obtain, under the hypothesis of Theorem \ref{t2}, the scaling limit of the process $X_t^n$.

\begin{theorem}
Assume that the hypothesis of Theorem \ref{t2} are fulfilled. Assume in addition that $d=1$, $\alpha>1$ and that 
\begin{itemize}
\item[{\bf (SG)}] There are positive constants $\kappa_0$, $k_0$ such that $|g(n+k_0)-g(n)| \geq \kappa_0$ for any $n \in \bb N$.
\end{itemize}

Then the tagged particle process satisfies
\[
\lim_{n \to \infty} \frac{X_t^n}{n} = \mc Z_t
\]
in distribution with respect to the $J$-Skorohod topology of $\mc D([0,\infty),\bb R)$, where $\mc Z_t$ is the unique process such that
\[
\exp\Bigg\{ i \theta \mc Z_t - c|\theta|^\alpha \int_0^T \frac{\phi(u(s,\mc Z_s))}{u(s,\mc Z_s)}ds\Bigg\}
\]
is a martingale for any $\theta \in \bb R$.
\end{theorem}

The constant $c$ in this theorem is equal to $\int h(x)\{e^{ix}-1\}dx$. Since we are restricted to dimension $d=1$, we can assume $h(x) = 1/|x|^{1+\alpha}$. Observe that in this theorem a superdiffusive scaling $t \to tn^\alpha$ is already embedded in the definition of $X_t^n$. The proof of this theorem follows from a careful adaptation of the arguments in \cite{Jar}, \cite{JLS}. All the needed tools we already introduced in the previous sections. We give here an sketch of proof, leaving the details to the reader.

\begin{proof}[Sketch of proof]

Instead of considering the pair $(\xi_t^n,X_t^n)$, we will define an auxiliary process $\zeta_t^n$ by taking $\zeta_t^n(x) = \xi_t^n(x+X_t^n)$ for $x \neq 0$ and $\zeta_t^n(0) = \xi_t^n(X_t^n)-1$. The process is known in the literature as the {\em environment as seen by the tagged particle} and it was introduced by Kipnis and Varadhan in \cite{KV}. Notice that the evolution of $\zeta_t^n$ is Markovian. Moreover, the {\em Palm measures} $\bar \nu_\rho(d\xi) = \rho^{-1}\xi(0)\nu_\rho(d\xi)$ are invariant under the evolution of $\zeta_t^n$. The evolution of $\zeta_t^n$ can be described as follows. At each site $x \neq 0$, a particle leaves $x$ and goes to $x+y$ at rate $n^\alpha p(y) g(\zeta_t^n(x))$. At site $x=0$, a particle leaves to site $y$ at rate $n^\alpha p(y) \bar g(\zeta_t^n(0))$, where $\bar g(n) = n g(n+1)/(n+1)$. And finally the whole configuration of particles is translated by $z$ at rate $n^\alpha p(z) b(\zeta_t^n(0))$, where $b(n) =g(n)/n$. The translations of the system correspond to jumps of the tagged particle. Let $N_t^{z,n}$ be the number of translations by $z$ performed by $\zeta_t^n$ up to time $t$. The number $N_t^{z,n}$ corresponds to the number of jumps by $z$ of the tagged particle up to time $t$. We have the formula
\[
\frac{X_t^n}{n} =\frac{1}{n} \sum_{z \in \bb Z_*} z N_t^{z,n},
\]
where $\bb Z_* = \bb Z \setminus \{0\}$.

Since $N_t^{z,n}$ is a Poisson process, the process 
\[
N_t^{z,n} -\int_0^t n^\alpha p(z) b(\zeta_t^n(0))ds 
\]
is a martingale. Notice that the compensators for $zN_t^{z,n}$ and $-zN_t^{-z,n}$ cancel each other. After checking the integrability of $X_t^n$, we conclude that $X_t^n$ is a martingale. More convenient will be to consider some exponential martingales associated to $X_t^n$. For each $\theta \in \bb R$ define
\begin{equation}
\label{ec30}
\mc M_t^{\theta,n} =
	\exp\Big\{ i\theta X_t^n/n -n^\alpha \sum_{z \in \bb Z_*} p(z) \big(e^{i\theta z /n} -1\big) \int_0^t b(\zeta_s^n(0))ds\Big\}.
\end{equation}

Notice that the integrand does not depend on $z$. The sum over $z$ of the expression before the integral is equal to
\[
\frac{1}{n} \sum_{z \in \bb Z_*} h(z/n) \{e^{i\theta z/n}-1\},
\]
which is a Riemann sum for $\int h(x) \{ e^{i\theta x}-1\}dx = -c|\theta|^\alpha$. The proof goes as follows. Define the empirical measure $\hat \pi_t^n$ by
\[
\hat \pi_t^n(dx) = \frac{1}{n} \sum_{z \in \bb Z} \zeta_t^n(z) \delta_{z/n}(dx).
\]

As in \cite{Jar}, the sequences $\{\hat \pi_\cdot^n\}_n$, $\{X_\cdot^n/n\}_n$ are tight. 
Notice that, aside of a factor $1/n$ corresponding to the tagged particle, the empirical measure $\hat \pi_t^n$ is the translation by $X_t^n/n$ of the empirical measure $\pi_t^n$ introduced before. 
We can find a subsequence $n'$ such that these two processes converge to some limits. We already know that $\pi_t^n(dx)$ converges to $u(x,t)dx$. 
Since the limit $u(x,t)$ is deterministic, we have joint convergence of the couple $\{(\pi_t^{n'},X_t^{n'}/n')\}_{n'}$ to $\{(u(x,t)dx, \mc Z_t)\}$ for some process $\mc Z_t$. We need to characterize the limiting process $\mc Z_t$. Observe that $\hat \pi_t^{n'}(dx)$ converges to $u(x+\mc Z_t,t)dx$. 
Notice that the expectation of $b(\zeta(0))$ with respect to $\bar \nu_\rho$ is equal to $\beta(\rho)=:\phi(\rho)/\rho$. In \eqref{ec30}, we want to replace $b(\zeta_t^n(0))$ by $\phi(\zeta_t^{n,\epsilon n}(0))/\zeta_t^{n,\epsilon n}(0)$, where $\zeta_t^{n,\epsilon n}(0) = (\epsilon n)^{-1} \sum_{k=1}^{\epsilon n} \zeta_s^n(k)$, since $\zeta_t^{n,\epsilon n}(0)$ is a function of the empirical measure $\hat \pi_t^n$. Such substitution is known as a {\em local replacement}, since it does not involve averaging with respect to a test function. The proof of this replacement confines us to dimension $d=1$ and $\alpha >1$ (see \cite{JLS} for a detailed discussion). In the proof of this local replacement, an spectral gap estimate, uniform on the density, for the process restricted to finite boxes is needed. At this point is where {\bf (SG)} is needed \cite{LSV}. The only difference with respect to \cite{JLS} is in the proof of the two-blocks estimate, where the corresponding local version of the moving particle lemma of Section \ref{s5} needs to be invoked. 

After proving this local replacement, we can pass to the limit in \eqref{ec30} to show that $\mc M_t^{\theta,n'}$ converges to 
\[
\mc M_t^\theta = \exp\Big\{ i\theta \mc Z_t + c|\theta|^\alpha \int_0^t \beta(u(s,\mc Z_s))ds\Big\}.
\]

Notice that $\beta(\rho) \leq \kappa$ for any $\rho \geq 0$. From this bound we can deduce that the sequence of martingales $\{\mc M_t^{\theta,n'}\}_{n'}$ is uniformly integrable, from where we obtain that $\mc M_t^\theta$ is a martingale. 

\end{proof}

\section*{Acknowledgements}

M.J. would like to thank Luis Silvestre and Cyril Imbert for stimulating discussions about fractional Laplacians.

\appendix

\section{The case $\alpha =2$}
\label{A2}

In this Appendix we explain how to obtain the hydrodynamic limit for the models with long jumps in the case $\alpha=2$. For simplicity, we restrict ourselves to dimension $d=1$; the arguments are the same in any dimension. Let us consider the transition rate $p(y-x) = c/|y-x|^3$. The right scaling for the exclusion process or the zero-range process associated to $p(\cdot)$ is $n^2/\log n$ (still superdiffusive). With this scaling, it turns out that the operator $\mc L_n$ defined in (\ref{ec4.0}) satisfies
\begin{equation}
\label{ec26}
\lim_{n \to \infty} \sup\big|L_n G(x/n) - c G''(x/n)\big| = 0,
\end{equation}

In fact, let $\epsilon >0$ be fixed. Remember that
\[
\mc L_n G(x/n) = \frac{1}{2n \log n} \sum_{z \in \bb Z^d} c\big| z/n\big|^{-3} \Big\{ G\Big(\frac{x+z}{n}\Big) + G\Big(\frac{x-z}{n}\Big)- 2G\Big(\frac{x}{n}\Big)\Big\}.
\]

Let us split the sum in two parts: when $|z| \leq \epsilon n$ and when $|z| \geq \epsilon n$. When $|z| \geq \epsilon n$, we estimate the term involving $G$ by $4\| G\|_\infty$, and we obtain the bound
\begin{multline*}
\frac{1}{n \log n} \sum_{|z| \geq \epsilon n} c\big| z/n\big|^{-3} \Big\{ G\Big(\frac{x+z}{n}\Big) + G\Big(\frac{x-z}{n}\Big)- 2G\Big(\frac{x}{n}\Big)\Big\}
	\leq \\ 
	\leq \frac{C(G)n^2}{\log n} \sum_{|z| \geq \epsilon n} c|z|^{-3} \leq \frac{C(G)}{\epsilon^2 \log n}.
\end{multline*}

When $|z| \leq \epsilon n$, we use a Taylor expansion to write the term involving $G$ as $G''(x/n) z^2/2n^2 + R^n_z z^4/n^4$, where $R^n_z$ is uniformly bounded. We conclude that
\begin{multline*}
\Big|\frac{1}{n \log n} \sum_{|z| \leq \epsilon n} \big| z/n\big|^{-3}  \Big\{ G\Big(\frac{x+z}{n}\Big) + G\Big(\frac{x-z}{n}\Big)- 2G\Big(\frac{x}{n}\Big)\Big\} - G''(x/n)\Big| \leq \\ 
	\leq \Big| G''(x/n)\Big(1-\frac{1}{\log n} \sum_{z=1}^{\epsilon n} \frac{1}{z}\Big)\Big|
\end{multline*}
plus a rest bounded by $C(G) \epsilon^2/\log n$. 
Since $\log(N+1) \leq \sum_{z=1}^N 1/z \leq \log(N+1) +1$, the right-hand side of the previous inequality goes to $0$ when $n \to \infty$. In order to prove the hydrodynamic limits of the processes $\eta_t^n$ and $\xi_t^n$, only suitable properties for the approximation operators $\mc L_n$ are needed, like the one we just proved. It is not difficult to show that $\mc L_n$ fulfills all the properties required.

\section{Condition {\bf (C)} and unbounded initial profiles}
\label{B}
In this Section we explain how to get rid of condition {\bf (C)} when $g(\cdot)$ is non-decreasing. Remember that for $\rho \leq \sigma$ we have  $\nu_\rho \preceq \nu_\sigma$. This is also true for the one-site marginals $q_\rho$, $q_\sigma$ as proved in Section \ref{s1.4}. Let us denote by $q_{\rho,\sigma}$ the coupling between $q_\rho$ and $q_\sigma$ constructed in Section \ref{s1.4}. In other words, $q_{\rho,\sigma}$ is a probability measure in $\bb N_0 \times \bb N_0$ with marginals $q_\rho$, $q_\sigma$ and satisfying $q_{\rho,\sigma}(x^1 \leq x^2)=1$. For simplicity, we will concentrate ourselves in the product initial measures described in Section \ref{s1.5}. The idea is to construct a more refined version of the coupling between two copies of the zero-range process described in Section \ref{s1.4}. 

Taking a closer look at the coupling $(\xi_t^1,\xi_t^2)$ constructed in Section \ref{s1.4}, we say that $\xi^1_t(x)$ is the number of {\em first-class particles} at site $x$ at time $t$. We say that $\xi_t^2(x)-\xi_t^1(x)$ is the number of {\em second-class particles} at site $x$ at time $t$. The idea is that first class particles do not feel the presence of second-class particles, while the evolution of second-class particles is modified by the presence of first-class particles at the same site. We will consider 4 types of particles: blue, green, red and white particles. The dynamics is as follows. We will start with a configuration with no green particles, and such that red and white particles do not share a site. Blue particles are first-class particles, green particles are second-class particles and red and white particles are third-class particles. Each time a red particle jumps over a site with at least one white particle, the red particle and one white particle are annihilated and a green particle is created. The same happens if a white particle jumps over a red particle. 

Let us define this process in a more precise way. We denote by $\xi_t^B(x)$, $\xi_t^G(x)$, $\xi_t^R(x)$ and $\xi_t^W(x)$ respectively the number of blue, green, red and white particles at site $x$ at time $t$. The process $(\xi^B_t,\xi^G_t,\xi^R_t,\xi^W_t)$ is defined in a proper subset of $\Omega_{zr}^4$ and it is generated by the operator $\bb L = \bb L^B +\bb L^G+\bb L^R+\hat{ \bb L}$, where
\begin{align*}
\bb L^B f(\xi^B,\xi^G,\xi^R,\xi^W)
	&= \sum_{x,y \in \bb Z^d} p(y-x) g(\xi^B(x))\times \\
	&\quad \quad \times\big\{f((\xi^B)^{x,y},\xi^G,\xi^R,\xi^W)-f(\xi^B,\xi^G,\xi^R,\xi^W)\big\},
\end{align*}
\begin{align*}
\bb L^G f(\xi^B,\xi^G,\xi^R,\xi^W)
	&=\sum_{x,y \in \bb Z^d} p(y-x)\big(g(\xi^G(x)+\xi^B(x))-g(\xi^B(x))\big)\times \\
	&\quad \quad \times \big\{f(\xi^B,(\xi^G)^{x,y},\xi^R,\xi^W)-f(\xi^B,\xi^G,\xi^R,\xi^W)\big\},
\end{align*}
\begin{align*}
\bb L^R f(\xi^B,\xi^G,\xi^R,\xi^W)
	&=\sum_{x,y \in \bb Z^d} p(y-x)\mathbf{1}(\xi^W(y)=0)\times \\
	&\quad \quad \times \big(g(\xi^G(x)+\xi^B(x)+\xi^R(x))-g(\xi^B(x)+\xi^G(x))\big)\times \\
	&\quad \quad \times \big\{f(\xi^B,\xi^G,(\xi^R)^{x,y},\xi^W)-f(\xi^B,\xi^G,\xi^R,\xi^W)\big\},
\end{align*}
\begin{align*}
\bb L^W f(\xi^B,\xi^G,\xi^R,\xi^W)
	&= \sum_{x,y \in \bb Z^d} p(y-x)\mathbf{1}(\xi^R(y)=0)\times \\
	&\quad \quad \times \big(g(\xi^G(x)+\xi^B(x)+\xi^W(x))-g(\xi^B(x)+\xi^G(x))\big)\times \\
	&\quad \quad \times \big\{f(\xi^B,\xi^G,\xi^R,(\xi^W)^{x,y})-f(\xi^B,\xi^G,\xi^R,\xi^W)\big\},
\end{align*}
\begin{align*}
\hat{\bb L} f(\xi^B,\xi^G,\xi^R,\xi^W)
	&=\sum_{x,y \in \bb Z^d} p(y-x)\mathbf{1}(\xi^W(y)\neq 0) \times \\
	&\quad \quad \times \big(g(\xi^G(x)+\xi^B(x)+\xi^R(x))-g(\xi^B(x)+\xi^G(x))\big)\times \\
	&\quad \quad \times \big\{f(\xi^B,\xi^G+\delta_y,\xi^R-\delta_x,\xi^W-\delta_y)-f(\xi^B,\xi^G,\xi^R,\xi^W)\big\}\\
	&+\sum_{x,y \in \bb Z^d} p(y-x)\mathbf{1}(\xi^R(y)\neq 0) \times \\
	&\quad \quad \times \big(g(\xi^G(x)+\xi^B(x)+\xi^W(x))-g(\xi^B(x)+\xi^G(x))\big)\times \\
	&\quad \quad \times \big\{f(\xi^B,\xi^G+\delta_y,\xi^R-\delta_y,\xi^W-\delta_x)-f(\xi^B,\xi^G,\xi^R,\xi^W)\big\}.
\end{align*}

The first line governes the evolution of blue particles and it is independent of the other particles. The second line governes the evolution of green particles, and it is independent of the evolution of red and white particles, except for the creation of new green particles, which corresponds to the last two lines. The third and fourth lines correspond to the evolution of white and red particles when they do not meet. 

Let $u_0(\cdot)$ be an unbounded profile, satisfying the hypothesis of Theorem \ref{t2}. Let us recall the definition of $\mu^n$: the measure $\mu^n$ is the product measure in $\Omega_{zr}$ of marginals given by $\mu^n(\xi(x)=l)= q_{u_0^n(x)}(l)$. Fix $M>0$. For each $n$, let us consider the product measure $\bar \mu^n$ in $\Omega_{zr}^2$ defined by
\[
\bar \mu^n(\xi^1(x) =l_1, \xi^2(x)=l_2) = \nu_{u_x^{M,n},U_x^{M,n}}(l_1,l_2),
\]
where $u_x^{M,n} = \min\{u_0^n(x),M\}$ and $U_x^{M,n} = \max\{u_0^n(x),M\}$. Now we define the initial measure $\hat \mu^n$ in $\Omega_{zr}^4$ by taking
\[
\xi^B_0(x) = \xi^1(x), \xi^G(x)=0,
\]
\[
\xi^R_0(x) =
\begin{cases}
0, & \text{Êif } u_0(x/n) \leq M\\
\xi^2(x) -\xi^1(x), &\text{Êif }Êu_0(x/n) > M,\\
\end{cases}
\]
\[
\xi^W_0(x) =
\begin{cases}
\xi^2(x) -\xi^1(x), & \text{Êif } u_0(x/n) \leq M\\
0, &\text{Êif }Êu_0(x/n) > M,\\
\end{cases}
\]

After a careful checking, we see that $\xi_t^B$ evolves like a zero-range process with initial distribution $\mu^{n,M}$ associated to the profile $\min\{u_0(x),M\}$, $\xi_t^B+\xi_t^G+\xi_t^R$ evolves like a zero-range process with initial distribution $\mu^n$ and that $\xi_t^B+\xi_t^G+\xi_t^W$ evolves like a zero-range process with initial distribution $\nu_M$. Therefore, Theorem \ref{t2} applies for the rescaled process $\xi_{tn^\alpha}^B$. Sending $M$ to infinity and using the monotone convergence theorem, we see that the rescaled process $\xi_{tn^\alpha}$ with initial condition $\mu^n$ also satisfies the conclusion of Theorem \ref{t2}.

\section{Condition {\bf (H)} and bounded profiles}
\label{C}
In this section we explain how to get rid of condition {\bf (H)} when $g(\cdot)$ is non-decreasing and the initial profile is bounded. Let $u_0(\cdot)$ be a bounded profile and let $\{\mu^n\}_n$ be the sequence of measures defined in Section \ref{s1.4}. Assume that there are constants $0<\rho_0<\rho_1$ such that $\rho_0 \leq u_0(x) \leq \rho_1$ for any $x \in \bb R^d$. We will use a coupling similar to the one used in the previous section. The idea is the following. Take $M >0$ and define 
\[
u_0^{M,+}(x) = 
\begin{cases}
u_0(x),& |x| \leq M\\
\rho_1, & |x| > M,
\end{cases}
\]
\[
u_0^{M,-}(x) = 
\begin{cases}
u_0(x),& |x| \leq M\\
\rho_0, & |x| > M,
\end{cases}
\]

We have $u_0^{M,-}(x) \leq u_0(x) \leq u_0^{M,+}(x)$ for any $x \in \bb R^d$. We will have three classes of particles: blue, green and red particles. The initial configurations are as follows. For any three numbers $\rho_0 \leq \rho \leq \rho_1$ it is possible to construct a measure $q_{\rho_0,\rho,\rho_1}$ in $\bb N \times \bb N \times \bb N$ with marginals $q_{\rho_0}$, $q_\rho$, $q_{\rho_1}$ and such that $q(x^1 \leq x^2 \leq x^3)=1$. For each $n$, let us consider the product measures $\tilde \mu^n$ in $\Omega_{zr}^3$ defined by
\begin{multline*}
\tilde \mu^n(\xi^B(x)=l_B, \xi^G(x) = l_G, \xi^R(x) =l_R) =\\
 =q_{u_0^{M-}(x/n), u_0(x/n),u_0^{M,+}(x/n)}(l_B,l_B+l_G,l_B+l_G+l_R).
\end{multline*}
The process $(\xi_t^B,\xi_t^G,\xi_t^R)$ with initial distribution $\tilde \mu^n$ is generated by the operator $\bb L$ defined in the previous section, projected on the set of configurations without white particles. As we pointed out before, $\xi_t^B$ evolves like a zero-range process with initial distribution associated to the profile $u_0^{M,-}$, $\xi_t^B+\xi_t^G$ evolves like a zero-range process with initial distribution $u_0$ and $\xi_t^B+\xi_t^G+\xi_t^R$ evolves like a zero-range process with initial distribution $u_0^{M,+}$. Both rescaled processes $\xi_{tn^\alpha}^B$ and $\xi_{tn^\alpha}^B+\xi_{tn^\alpha}^G+\xi_{tn^\alpha}^R$ satisfy the hypothesis of Theorem \ref{t2}. Tightness of the empirical measure associated to $\xi_t^R+\xi_t^G$ follows by comparison, and the limiting measures $\pi_t(dx)=u(x,t)dx$ (possibly random) are bounded between $u^{M,-}(x,t)dx$ and $u^{M,+}(x,t) dx$, which are the solutions of the hydrodynamic equation \eqref{ec2} with initial conditions $u_0^{M,-}$, $u_0^{M,+}$. Sending $M$ to $\infty$ we obtain that the rescaled process $\xi_{tn^\alpha}^B+\xi_{tn^\alpha}^G$ also satisfies the conclusion of Theorem \ref{t3}. When $u_0$ is not bounded below, a similar, but more sophisticated coupling can be constructed.



\end{document}